\documentclass[10pt,sort]{elsarticle}
\usepackage[a4paper]{geometry}
\usepackage{graphicx}
\usepackage{amssymb}
\usepackage{epstopdf}
\usepackage{color}
\usepackage{float}
\usepackage{stmaryrd}
\usepackage{mathrsfs}
\usepackage{amsmath}
\usepackage{booktabs}
\usepackage{amsthm}
\usepackage{bm}
\usepackage{setspace}
\usepackage{cancel}
\theoremstyle{plain}

\theoremstyle{definition}

\theoremstyle{remark}
\newtheorem{rem}{Remark}

\newcommand{\jump}[1]{\left\llbracket#1\right\rrbracket}

\numberwithin{equation}{section}

\begin{document}
\begin{frontmatter}

\title{A Provably Stable Discontinuous Galerkin Spectral Element Approximation for Moving Hexahedral Meshes}
\author[fsu]{David A. Kopriva}
\author[unikoeln]{Andrew R. Winters\corref{cor1}}
\cortext[cor1]{Corresponding author.}
\ead{awinters@math.uni-koeln.de}
\author[unikoeln]{Marvin Bohm}
\author[unikoeln]{Gregor. J. Gassner}

\address[fsu]{Department of Mathematics, The Florida State University, Tallahassee, FL 32306, USA}
\address[unikoeln]{Mathematisches Institut, Universit\"{a}t zu K\"{o}ln, Weyertal 86-90, 50931 K\"{o}ln, Germany}

\begin{abstract}
We design a novel provably stable discontinuous Galerkin spectral element (DGSEM) approximation to solve systems of conservation laws on moving domains. To incorporate the motion of the domain, we use an arbitrary Lagrangian-Eulerian formulation to map the governing equations to a fixed reference domain. The approximation is made stable by a discretization of a skew-symmetric formulation of the problem. We prove that the discrete approximation is stable, conservative and, for constant coefficient problems, maintains the free-stream preservation property. We also provide details on how to add the new skew-symmetric ALE approximation to an existing discontinuous Galerkin spectral element code. Lastly, we provide numerical support of the theoretical results. 
\end{abstract}


\end{frontmatter}
\section{Introduction}

Many applications in computational physics require the numerical approximation of a system of hyperbolic conservation laws on moving domains, for example problems in fluid dynamics \cite{etienne2009perspective,Lomtevetal1999,Mavriplis:2006mb,wang2015high} or electromagnetism \cite{censor2004,harfoush1989,ho2006}. A common approach to approximate solutions of partial differential equations (PDEs) with moving boundaries is to use an arbitrary Lagrangian-Eulerian (ALE) formulation \cite{etienne2009perspective,Hay2014204}. The ALE method maps a time dependent domain, $\Omega$, enclosed by the moving boundaries onto a fixed reference domain, $E$. Conveniently, systems of conservation laws in the original moving domain are transformed to a set of conservation law equations in the reference domain. In the numerical approximation on the reference domain, the new set of equations depends on the mesh velocity.

A discontinuous Galerkin spectral element method for moving domains (DGSEM-ALE) that was spectrally accurate in space, free-stream preserving, and retained the full time accuracy of the time integrator was introduced in \cite{Minoli:2010rt}. Extensions of the method were presented in \cite{ISI:000302933600007}, and \cite{Winters:2013nx,Winters:2013oq,Winters:2014hl}, the latter of which addressed the approximation of problems with discontinuous and moving material interfaces and efficiency through the addition of local time stepping.

None of the papers on the DGSEM-ALE approximation directly addressed stability, however. In fact, it has been noted that even on static domains, discontinuous Galerkin spectral element approximations for variable coefficient problems \cite{Beck:2013sf}, and problems that become variable coefficient because of curved elements \cite{kopriva2015}, can be unstable.


In this paper, we now address the additional issues of robustness and the ability of a moving mesh method to be guaranteed stable. Recent work for static meshes has focused on the use of DGSEM approximations written in a skew-symmetric form to guarantee stability, e.g. \cite{gassner2015,Gassner:2013uq,kopriva2015}. The skew-symmetric form of a problem is usually found by averaging the conservative form and a non-conservative advective form of the equation. This is problematic, because it is not obvious that discretizations of the skew-symmetric form remain conservative. Recent success has been had using diagonal norm summation-by-parts (SBP) operators to discretize the spatial derivatives in the skew-symmetric formulation \cite{skew_sbp2,gassner_skew_burgers,gassner_kepdg,gassner2015}. There is a known link between SBP methods and the discontinuous Galerkin spectral element approximation with Legendre-Gauss-Lobatto points, e.g. \cite{gassner_skew_burgers}. The approximation developed here will also use a skew-symmetric formulation.

In addition to stability, the DGSEM approximation that we propose is conservative, high order accurate in both space and time, and ensures that for constant coefficient problems a constant solution of the conservation law remains constant, discretely, for all time, i.e, the approximation possesses free-stream preservation. Failure to satisfy free-stream preservation usually implies that the motion of the mesh can create spurious waves and may introduce discrete errors that can lead to wave misidentification, even for spectrally accurate approximations \cite{Kopriva:2006er}. 

A stable deforming mesh approximation for high-order finite difference schemes with the summation by parts (SBP) property was recently proposed by Nikkar and Nordstr\"{o}m \cite{Nikkar:2014si}. The method development here parallels theirs, and the result satisfies the same type of energy estimate, even though our use of the weak rather than strong formulation, the notation, and the approximation differ.


The remainder of this paper is organized as follows:  Sec. \ref{ALESection} reviews the ALE mapping between a reference space and a general curvilinear coordinate system. Next, in Sec. \ref{sec:WellPosedness}, we demonstrate the well-posedness of the ALE formulation, which presents the target to be approximated. In Sec. \ref{StableFormulation} we use a skew-symmetric formulation of the governing equations to develop a stable approximation of the problem on moving meshes. Sec. \ref{PropertiesOfStable} provides proofs of the conservation, stability, and free-stream preserving properties of the skew-symmetric DGSEM-ALE formulation. We provide some details on the implementation of the newly proposed scheme in Sec. \ref{ImplementationSection}. Numerical results are presented in Sec. \ref{NumericalResultsSection} to support the theoretical findings. Finally, concluding remarks are given in Sec. \ref{ConclusionSection}.

\section{Arbitrary Lagrangian-Eulerian (ALE) Formulation of Conservation Laws in a Curvilinear Coordinate System}\label{ALESection}

We will derive a discontinuous Galerkin spectral element method for systems of partial differential equations of the form
\begin{equation}
{{{\mathbf{q}}}_t} + \nabla  \cdot \vec {\mathbf{f}} = {{{\mathbf{q}}}_t} + \sum\limits_{i = 1}^3 {\frac{{\partial {{ {\mathbf{f}}}_i}}}{{\partial {x_i}}} = 0},
\label{eq:OrigConsLaw}
\end{equation}
on a three dimensional domain with moving boundaries, $\Omega\left(\vec x,t\right)$, where $\vec x = \left(x_{1},x_{2},x_{3}\right)= (x,y,z)$. Here we denote the solution and flux vector components by bold face and spatial vectors by arrows.
We assume that the system is symmetric hyperbolic, with covariant flux components
\begin{equation}
 {\mathbf{f}}_{i} = A_{i}(\vec x){\mathbf{q}},
\end{equation}
where the $A_{i} = A^{T}_{i}$ are matrices. If the system is not symmetric to start with, then symmetrization, which is available for most systems of interest \cite{mccarthy1980,R.-F.-Warming:1975fk}, is applied first. We further assume that the matrices are smooth, having bounded derivatives. Since the system is hyperbolic, there exists an invertible matrix $P(A)$ such that 
\begin{equation}
A=\sum\limits_{i = 1}^3 {{\alpha _i}{A_i}}  = P(A)\Lambda(A) {P^{ - 1}(A)},\end{equation}
for any $\vec \alpha$ with $\sum\limits_{i = 1}^3 {\alpha _i^2}  \ne 0$ and some real diagonal matrix $\Lambda$.

As a concrete example, the symmetric wave equation can be written in the form (\ref{eq:OrigConsLaw}) as
\begin{equation}
{\left[ {\begin{array}{*{20}{c}}
  p \\ 
  u \\ 
  v \\ 
  w 
\end{array}} \right]_t} + {\left( {\left[ {\begin{array}{*{20}{c}}
  0&c&0&0 \\ 
  c&0&0&0 \\ 
  0&0&0&0 \\ 
  0&0&0&0 
\end{array}} \right]\left[ {\begin{array}{*{20}{c}}
  p \\ 
  u \\ 
  v \\ 
  w 
\end{array}} \right]} \right)_x} + {\left( {\left[ {\begin{array}{*{20}{c}}
  0&0&c&0 \\ 
  0&0&0&0 \\ 
  c&0&0&0 \\ 
  0&0&0&0 
\end{array}} \right]\left[ {\begin{array}{*{20}{c}}
  p \\ 
  u \\ 
  v \\ 
  w 
\end{array}} \right]} \right)_y} + {\left( {\left[ {\begin{array}{*{20}{c}}
  0&0&0&c \\ 
  0&0&0&0 \\ 
  0&0&0&0 \\ 
  c&0&0&0 
\end{array}} \right]\left[ {\begin{array}{*{20}{c}}
  p \\ 
  u \\ 
  v \\ 
  w 
\end{array}} \right]} \right)_z} = 0.
\end{equation}

In the ALE formulation, one maps $\Omega\left(\vec x,t\right)$ onto a reference domain, $E$, by a transformation
\begin{equation}
\left(\vec x, t\right) = \vec X\left(\vec \xi,\tau\right),
\label{eq:mapping}
\end{equation}
where $\vec \xi = \left(\xi^{1},\xi^{2},\xi^{3}\right) = \left(\xi,\eta,\zeta\right)$ is a three dimensional curvilinear coordinate on the reference domain. 
For convenience with the approximations later, we can take the reference domain to be the reference cube $E = [-1,1]^{3}$. The mapping has a set of covariant basis vectors, $\vec a_{i}$ defined 
by
\begin{equation}
{{\vec a}_i} = \frac{{\partial \vec X}}{{\partial {\xi ^i}}}\quad i = 1,2,3.
\end{equation}

From the covariant basis, one can formally define the contravariant basis $\vec a^{i}$, multiplied by the Jacobian of the transformation, $\mathcal{J}$,
 by
\begin{equation}
\mathcal{J}\vec a^i  = \mathcal{J}\nabla \xi ^i  = \vec a_j  \times \vec a_k  = \frac{{\partial \vec X}}
{{\partial \xi ^j }} \times \frac{{\partial \vec X}}
{{\partial \xi ^k }},\quad \left( {i,j,k} \right)\,cyclic.
\label{eq:contravaraiantBasis}
\end{equation}
The Jacobian itself can be written in terms of the covariant vectors as
\begin{equation}
	\mathcal{J} = \vec a_i  \cdot \left( {\vec a_j  \times \vec a_k } \right),\quad \left( {i,j,k} \right)\;cyclic.
\end{equation}
The geometry satisfies the well-known metric identities \cite{Vinokur1974}, \cite[Chpt. III]{Thompsonetal1985}
\begin{equation}
	\sum\limits_{i = 1}^3 {\frac{{\partial {\mathcal{J}}{{\vec a}^i}}}{{\partial {\xi ^i}}}}  = 0,
\label{eq:MetricIdentities}
\end{equation}
and the Geometric Conservation Law (GCL) \cite{Mavriplis:2006mb}
\begin{equation}
\frac{{\partial \mathcal{J}}}{{\partial \tau }} - \sum\limits_{i = 1}^3 {\frac{{\partial \mathcal{J}{{\vec a}^i} \cdot {{\vec x}_\tau }}}{{\partial {\xi ^i}}} = 0}. 
\label{eq:GCL}
\end{equation}

Under the transformation (\ref{eq:mapping}), see, e.g. \cite{Minoli:2010rt}, the conservation law (\ref{eq:OrigConsLaw})
remains a conservation law 
\begin{equation}
\frac{{\partial \mathcal{J}{\mathbf{q}}}}{{\partial \tau }} + \sum\limits_{i = 1}^3 {\frac{{\partial \mathcal{J}{{\vec a}^i} \cdot \left( {\vec {\mathbf{F}} - {\mathbf{q}}  {{\vec x}_\tau }} \right)}}{{\partial {\xi ^i}}} = 0}, 
\label{eq:transfconsLaw1}
\end{equation}
on the reference domain.
If we define the contravariant coefficient matrices, 
\begin{equation}
{{\tilde A}^i} = \mathcal{J}{{\vec a}^i} \cdot \sum\limits_{j = 1}^3 {{A_j}{{\hat x}_j}},\quad i = 1,2,3,
    \label{eq:ContravariantCoefMatrices}
\end{equation}
where $\hat x_{j}$ denotes the unit vector in the $j^{th}$ coordinate direction, then we can re-write (\ref{eq:transfconsLaw1}) in fully conservative form as
\begin{equation}
\frac{{\partial \mathcal{J}{\mathbf{q}}}}{{\partial \tau }} + \sum\limits_{i = 1}^3 {\frac{{\partial {\mathbf{\tilde f}}}}{{\partial {\xi ^i}}} = 0} ,
\end{equation}
where ${\mathbf{\tilde f}} = \left( {{{\tilde A}^i} - \mathcal{J}{{\vec a}^i} \cdot {{\vec x}_\tau }I} \right){\mathbf{q}}$ and $I$ is the identity matrix.

Using the metric identities, the GCL and the conservative form of the equations, we can also derive a non-conservative form of the equations. From the chain rule,
\begin{equation}{{\mathcal{J}}_\tau }{\mathbf{q}} + {\mathcal{J}}{{{\mathbf{q}}}_\tau } + \sum\limits_{i=1}^3 {\frac{{\partial {{\tilde A}^i}}}{{\partial {\xi ^i}}}{\mathbf{q}}}  + \sum\limits_{i=1}^3 {{{\tilde A}^i}\frac{{\partial {\mathbf{q}}}}{{\partial {\xi ^i}}}}  - \sum\limits_{i=1}^3 {\frac{{\partial {\mathcal{J}}{{\vec a}^i} \cdot {{\vec x}_\tau }}}{{\partial {\xi ^i}}}{\mathbf{q}}}  - \sum\limits_{i=1}^3 {{\mathcal{J}}{{\vec a}^i} \cdot {{\vec x}_\tau }\frac{{\partial {\mathbf{q}}}}{{\partial {\xi ^i}}}}  = 0.\end{equation}
Applying the GCL \eqref{eq:GCL} we find,
\begin{equation}
\mathcal{J}{{\mathbf{q}}_\tau } + \sum\limits_{i=1}^3 {{{\tilde A}^i}\frac{{\partial {\mathbf{q}}}}{{\partial {\xi ^i}}}}  - \sum\limits_{i=1}^3 {\mathcal{J}{{\vec a}^i} \cdot {{\vec x}_\tau }\frac{{\partial {\mathbf{q}}}}{{\partial {\xi ^i}}}}  + \sum\limits_{i=1}^3 {\frac{{\partial {{\tilde A}^i}}}{{\partial {\xi ^i}}}{\mathbf{q}}}  = 0,
\end{equation}
giving us the nonconservative system of equations
\begin{equation}{\mathcal{J}}{{{\mathbf{q}}}_\tau } + \sum\limits_{i=1}^3 {\left( {{{\tilde A}^i} - J{{\vec a}^i} \cdot {{\vec x}_\tau }I} \right)\frac{{\partial {\mathbf{q}}}}{{\partial {\xi ^i}}}}+ \sum\limits_{i=1}^3 {\frac{{\partial {{\tilde A}^i}}}{{\partial {\xi ^i}}}{\mathbf{q}}}  = 0.
\end{equation}

Finally, let us define the matrices
\begin{equation}{\mathcal{A}^i}\left( {\vec \xi ,\tau } \right) = {{\tilde A}^i} - J{{\vec a}^i} \cdot {{\vec x}_\tau }I,\quad i = 1,2,3,\end{equation}
to let us write the conservative and non-conservative forms of the equations as
\begin{equation}
\frac{{\partial \mathcal{J}{\mathbf{q}}}}{{\partial \tau }} + \sum\limits_{i = 1}^3 {\frac{\partial }{{\partial {\xi ^i}}}\left( {{\mathcal{A}^i}{\mathbf{q}}} \right) = 0}, 
\label{eq:ConsFormMapped}
\end{equation}
and
\begin{equation}\mathcal{J}{{{\mathbf{q}}}_\tau } + \sum\limits_{i=1}^3 {{\mathcal{A}^i}\frac{{\partial {\mathbf{q}}}}{{\partial {\xi ^i}}}}+ \sum\limits_{i=1}^3 {\frac{{\partial {{\tilde A}^i}}}{{\partial {\xi ^i}}}{\mathbf{q}}}  = 0,
\label{eq:NonConsFormMapped}
\end{equation}
respectively.

\section{Well-Posedness of the ALE Formulation}\label{sec:WellPosedness}

Since the stability proof mimics that of well-posedness of the PDE, 
we first show that the system of PDEs on the moving domain is well-posed when appropriate boundary conditions are applied. The derivation here follows that of \cite{Nikkar:2014si}, but uses notation that is consistent with how we will write the DGSEM and its stability proof.  

Let us define the inner 
product on the reference domain as
\begin{equation}
\left( {\mathbf{u},\mathbf{v}} \right) = \int_E {{{\mathbf{u}}^T}\mathbf{v}\,d\vec \xi } ,
\end{equation}
and norm $\left\| {\mathbf{u}} \right\| = \sqrt {\left( {\mathbf{u},\mathbf{u}} \right)}$. We denote the space $\mathbb{L}^{2} = \left\{\mathbf{u}: \left\| {\mathbf{u}} \right\|<\infty\right\}$.

We show well-posedness by bounding the time rate of change of the energy in the moving domain $\Omega$, which is equivalent to
\begin{equation}\frac{d}{{d\tau }}\left\| {{\mathbf{q}}} \right\|_\mathcal{J}^2 \equiv \frac{d}{{d\tau }}\int\limits_E {{{{\mathbf{q}}}^T}{\mathbf{q}}\mathcal{J}\,d\vec \xi }  = \left( {{\mathbf{q}},\frac{{\partial \mathcal{J}{\mathbf{q}}}}{{\partial \tau }}} \right) + \left( {{\mathbf{q}},\mathcal{J}{{{\mathbf{q}}}_\tau }} \right).
\label{eq:timeDerivOfenergy}
\end{equation}
The two terms on the right of (\ref{eq:timeDerivOfenergy}) can be found from the conservative and 
non-conservative forms of the mapped system.
The inner product of the conservative form of the equation (\ref{eq:ConsFormMapped}) with the solution is
\begin{equation}\label{eq:IPcons}\left( {{\mathbf{q}},\frac{{\partial \mathcal{J}{\mathbf{q}}}}{{\partial \tau }}} \right) + \sum\limits_{i = 1}^3 {\left( {{\mathbf{q}},\frac{\partial \left( {{\mathcal{A}^i}{\mathbf{q}}} \right)}{{\partial {\xi ^i}}}} \right) = 0}. \end{equation}
Similarly, using the nonconservative form \eqref{eq:NonConsFormMapped},
\begin{equation}\label{eq:IPNonCons}
\left( {{\mathbf{q}},J{{\mathbf{q}}_\tau }} \right) + \sum\limits_{i=1}^3 {\left( {{\mathbf{q}},{\mathcal{A}^i}\frac{{\partial {\mathbf{q}}}}{{\partial {\xi ^i}}}} \right)}  + \sum\limits_{i=1}^3 {\left( {{\mathbf{q}},\frac{{\partial {{\tilde A}^i}}}{{\partial {\xi ^i}}}{\mathbf{q}}} \right)}  = 0.
\end{equation}
Adding the conservative \eqref{eq:IPcons} and non-conservative \eqref{eq:IPNonCons} forms together gives
\begin{equation}
\left( {{\mathbf{q}},\frac{{\partial \mathcal{J}{\mathbf{q}}}}{{\partial \tau }}} \right) + \left( {{\mathbf{q}},J{{{\mathbf{q}}}_\tau }} \right) + \sum\limits_{i=1}^3 {\left( {{\mathbf{q}},\frac{\partial \left( {{\mathcal{A}^i}{\mathbf{q}}} \right)}{{\partial {\xi ^i}}} + {\mathcal{A}^i}\frac{{\partial {\mathbf{q}}}}{{\partial {\xi ^i}}}} \right)} + \sum\limits_{i=1}^3 {\left( {{\mathbf{q}},\frac{{\partial {{\tilde A}^i}}}{{\partial {\xi ^i}}}{\mathbf{q}}} \right)} = 0.
\label{eq:summedWeakForms}
\end{equation}
Also, because the matrices $\mathcal{A}^{i}$ are symmetric, we see that
\begin{equation}
\frac{\partial }{{\partial {\xi ^i}}}\left( {{{{\mathbf{q}}}^T}{\mathcal{A}^i}{\mathbf{q}}} \right) = \frac{\partial }{{\partial {\xi ^i}}}\left( {{{{\mathbf{q}}}^T}} \right){\mathcal{A}^i}{\mathbf{q}} + {{{\mathbf{q}}}^T}\frac{\partial \left( {{\mathcal{A}^i}{\mathbf{q}}} \right)}{{\partial {\xi ^i}}} = {{{\mathbf{q}}}^T}{\mathcal{A}^i}\frac{\partial }{{\partial {\xi ^i}}}\left( {{\mathbf{q}}} \right) + {{{\mathbf{q}}}^T}\frac{\partial \left( {{\mathcal{A}^i}{\mathbf{q}}} \right)}{{\partial {\xi ^i}}},
\label{eq:VolDivergence}
\end{equation}
so
\begin{equation}
	\sum\limits_{i=1}^3 {\left( {{\mathbf{q}},\frac{\partial }{{\partial {\xi ^i}}}\left( {{\mathcal{A}^i}{\mathbf{q}}} \right) + {\mathcal{A}^i}\frac{{\partial {\mathbf{q}}}}{{\partial {\xi ^i}}}} \right)}  = \int\limits_E {\sum\limits_{i=1}^3 {\frac{\partial }{{\partial {\xi ^i}}}\left( {{{{\mathbf{q}}}^T}{\mathcal{A}^i}{\mathbf{q}}} \right)}d\vec\xi }. \end{equation}
Then with (\ref{eq:summedWeakForms}) and (\ref{eq:VolDivergence}), (\ref{eq:timeDerivOfenergy}) can be written with the divergence of a flux as
\begin{equation}
\frac{d}{{d\tau }}\left\| {{\mathbf{q}}} \right\|_\mathcal{J}^2 + \int\limits_E {\sum\limits_{i=1}^3 {\frac{\partial }{{\partial {\xi ^i}}}\left( {{{{\mathbf{q}}}^T}{\mathcal{A}^i}{\mathbf{q}}} \right)\,d\vec \xi } } + \sum\limits_{i=1}^3 {\left( {{\mathbf{q}},\frac{{\partial {{\tilde A}^i}}}{{\partial {\xi ^i}}}{\mathbf{q}}} \right)} = 0.
\label{eq:EnergyEqn1}
\end{equation}
Now define
\begin{equation}\vec {\mathcal{A}}  = \sum\limits_{i=1}^3 {{\mathcal{A}^i}{{\hat \xi }^i}}, \end{equation}
where $\hat{\xi}^i$ are the cardinal bases and $\hat{n}$ the outward normal on the reference box. Then Gauss' theorem states that
\begin{equation}
\frac{d}{{d\tau }}\left\| {{\mathbf{q}}} \right\|_\mathcal{J}^2 + \int\limits_{\partial E } {{{{\mathbf{q}}}^T}\vec {\mathcal{A}} \cdot \hat n{\mathbf{q}}dS}=- \sum\limits_{i=1}^3 {\left( {{\mathbf{q}},\frac{{\partial {{\tilde A}^i}}}{{\partial {\xi ^i}}}{\mathbf{q}}} \right)}.
\end{equation}
Next, we determine the bound
\begin{equation} - \sum\limits_{i=1}^3 {\left( {{\mathbf{q}},\frac{{\partial {{\tilde A}^i}}}{{\partial {\xi ^i}}}{\mathbf{q}}} \right)}  = \left( {{\mathbf{q}},\left\{ { - \frac{1}{\mathcal{J}}\sum\limits_{i=1}^3 {\frac{{\partial {{\tilde A}^i}}}{{\partial {\xi ^i}}}} } \right\}\mathcal{J}{\mathbf{q}}} \right) \leqslant \mathop {\max }\limits_E \left| {\frac{1}{\mathcal{J}}\sum\limits_{i=1}^3 {\frac{{\partial {{\tilde A}^i}}}{{\partial {\xi ^i}}}} } \right|\left\| {\mathbf{q}} \right\|_\mathcal{J}^2.\end{equation}
We then note that under the transformation rules,
\begin{equation}\frac{1}{\mathcal{J}}\sum\limits_{i=1}^3 {\frac{{\partial {{\tilde A}^i}}}{{\partial {\xi ^i}}}}  = \sum\limits_{i=1}^3 {\frac{{\partial {A_i}}}{{\partial {x_i}}}}, \end{equation} is assumed to be bounded. Therefore,
\begin{equation}\frac{d}{{d\tau }}\left\| {\mathbf{q}} \right\|_\mathcal{J}^2 + \int\limits_{\partial E} {{{\mathbf{q}}^T}\vec{ \mathcal{A}  }\cdot \hat n{\mathbf{q}}\,dS}  \leqslant   2\gamma \left\| {\mathbf{q}} \right\|_\mathcal{J}^2,
\label{eq:EnergyEqnSurf}
\end{equation}
where
\begin{equation}\gamma  = \frac{1}{2}\mathop {\max }\limits_E \left| {\frac{1}{\mathcal{J}}\sum\limits_{i=1}^3 {\frac{{\partial {{\tilde A}^i}}}{{\partial {\xi ^i}}}} } \right| < \infty.
\end{equation}

Integrating both sides of (\ref{eq:EnergyEqnSurf}) with respect to time over a time interval $[0,T]$ gives
\begin{equation}
\left\| {{\mathbf{q}}\left( T \right)} \right\|_\mathcal{J}^2 \leqslant {e^{2\gamma T}}\left\{ {\left\| {{\mathbf{q}}\left( 0 \right)} \right\|_\mathcal{J}^2 - \int_0^T {\int_{\partial E} {{e^{ - 2\gamma \tau }}{{\mathbf{q}}^T}\vec{ \mathcal{A}  }\cdot \hat n{\mathbf{q}}\,dSd\tau } } } \right\}.
\label{eq:EnergyEqnSurfIntegrated}
\end{equation}

We now apply boundary conditions to (\ref{eq:EnergyEqnSurfIntegrated}). First we note that  $\vec {\mathcal{A}}\cdot\hat n$ is diagonalizable. From the definition,
\begin{equation}{\mathcal{A}^i} = J{{\vec a}^i} \cdot \left( {\sum\limits_{j=1}^3 {{A_j}{{\hat x}_j}}  - {{\vec x}_\tau }I} \right),\end{equation}
so the normal flux matrices are given by
\begin{equation}\sum\limits_{i=1}^3 {{\mathcal{A}^i}{n^{i}}}  = \sum\limits_{j=1}^3 {{A_j}\left\{ {\sum\limits_{i=1}^3 {J{{\vec a}^i} \cdot {{\hat x}_j}{{\hat n}^{i}}} } \right\}}  - \left\{ {\sum\limits_{i=1}^3 {J{{\vec a}^i} \cdot {{\hat x}_\tau }{{\hat n}^{i}}} } \right\}I = \sum\limits_{j=1}^3 {{A_j}{\alpha _j}}  - \beta I.\end{equation}
Since $\sum\limits_{j=1}^3 {{A_j}{\alpha _j}}$ is diagonalizable, so is $S\equiv\sum\limits_{i=1}^3 {{\mathcal{A}^i}{\hat n^{i}}} = P\left(S\right)\Lambda(S) P^{-1}(S)$. This means that we can split the boundary contributions into incoming and outgoing components according to the sign of the eigenvalues. So let us split the normal coefficient matrix
\begin{equation}\sum\limits_{i=1}^3 {{\mathcal{A}^i}{n^{i}}}  = P{\Lambda ^ + }{P^{ - 1}} + P{\Lambda ^ - }{P^{ - 1}} = A^{+} + A^{-},\end{equation}
so that
\begin{equation}
\left\| {{\mathbf{q}}\left( T \right)} \right\|_\mathcal{J}^2 \le {e^{2\gamma T}}\left\{ {\left\| {{\mathbf{q}}\left( 0 \right)} \right\|_\mathcal{J}^2 - \int_0^T {\int_{\partial E} {{e^{ - 2\gamma \tau }}{{\mathbf{q}}^T}P{\Lambda ^ + }{P^{ - 1}}{\mathbf{q}}\,dSd\tau } }  - \int_0^T {\int_{\partial E} {{e^{ - 2\gamma \tau }}{{\mathbf{q}}^T}P{\Lambda ^ - }{P^{ - 1}}{\mathbf{q}}\,dSd\tau } } } \right\}.
 \end{equation}
We then replace the incoming values (associated with the $\Lambda^{-}$ eigenvalues) with the exterior state ${\mathbf{q}}_{\infty}$  and bound the integrals over time to get
\begin{equation}
\left\| {{\mathbf{q}}\left( T \right)} \right\|_\mathcal{J}^2 + {e^{2\gamma T}}\int_0^T {\int_{\partial E} {{{\mathbf{q}}^T}P{\Lambda ^ + }{P^{ - 1}}{\mathbf{q}}\,dSd\tau } }  \leqslant {e^{2\gamma T}}\left\{\left\| {{\mathbf{q}}\left( 0 \right)} \right\|_\mathcal{J}^2 + \int_0^T {\int_{\partial E} {{\mathbf{q}}_\infty ^TP\left| {{\Lambda ^ - }} \right|{P^{ - 1}}{{\mathbf{q}}_\infty }\,dSd\tau } }\right\} .
\label{eq:StronglyWellPosed}
\end{equation}
The statement (\ref{eq:StronglyWellPosed}) says that the initial boundary value problem is \emph{strongly well posed} \cite{Lorenz:1989fk}. When we set the boundary states to zero, we see that the system of equations is \emph{well-posed},
\begin{equation}\left\| {{\mathbf{q}}\left( T \right)} \right\|_\mathcal{J} \leqslant {e^{\gamma T}}\left\| {{\mathbf{q}}\left( 0 \right)} \right\|_\mathcal{J}.\end{equation} Finally, if the matrices are also constant, then $\gamma = 0$ and the energy never grows,
\begin{equation}\left\| {{\mathbf{q}}\left( T \right)} \right\|_\mathcal{J} \leqslant \left\| {{\mathbf{q}}\left( 0 \right)} \right\|_\mathcal{J}.\end{equation}

\section{A Stable DGSEM-ALE for Moving Domains}\label{StableFormulation}

We now derive a discontinuous Galerkin spectral element method (DGSEM) for moving elements whose stability properties mimic (\ref{eq:StronglyWellPosed}).
A description of the standard approximation can be found in \cite{Kopriva:2009nx}. We subdivide the domain $\Omega$ into non-overlapping, geometrically conforming hexahedral elements $e^{m}$ that cover $\Omega$. Since $\Omega$ has moving boundaries, so too will the elements. We then map each element individually with a local time dependent mapping of the form (\ref{eq:mapping}) onto the reference element $E$. Then on each element, the equations take on the conservative form (\ref{eq:ConsFormMapped}) and the non-conservative form
 \begin{equation}
\frac{{\partial \mathcal{J}}}{{\partial \tau }}{\mathbf{q}} + \mathcal{J}\frac{{\partial {\mathbf{q}}}}{{\partial \tau }} + \sum\limits_{i = 1}^3 {\frac{{\partial {\mathcal{A}^i}}}{{\partial {\xi ^i}}}{\mathbf{q}}}  + \sum\limits_{i = 1}^3 {{\mathcal{A}^i}\frac{{\partial {\mathbf{q}}}}{{\partial {\xi ^i}}}}  = 0,
 \end{equation}
 which is written without applying the metric identities and GCL.
 
 To create the skew-symmetric form, we average the conservative and nonconservative forms \cite{Gassner:2013uq} to get
 \begin{equation}
\frac{1}{2}\left\{ {\frac{{\partial \mathcal{J}{\mathbf{q}}}}{{\partial \tau }} + \frac{{\partial \mathcal{J}}}{{\partial \tau }}{\mathbf{q}} + \mathcal{J}\frac{{\partial {\mathbf{q}}}}{{\partial \tau }}} \right\} + \frac{1}{2}\sum\limits_{i = 1}^3 {\frac{\partial\left( {{\mathcal{A}^i}{\mathbf{q}}} \right) }{{\partial {\xi ^i}}}}  + \frac{1}{2}\left\{ {\sum\limits_{i = 1}^3 {\frac{{\partial {\mathcal{A}^i}}}{{\partial {\xi ^i}}}{\mathbf{q}}}  + \sum\limits_{i = 1}^3 {{\mathcal{A}^i}\frac{{\partial {\mathbf{q}}}}{{\partial {\xi ^i}}}} } \right\} = 0.
\label{eq:Averagedform}
 \end{equation}
 We  construct a weak form of (\ref{eq:Averagedform}) by multiplying by a test function $\boldsymbol{\phi}\in\mathbb{L}^{2}$ and integrating over the domain. In inner product notation, the weak form is
 \begin{equation}
\frac{1}{2}\left( {\frac{{\partial \mathcal{J}{\mathbf{q}}}}{{\partial \tau }} + \frac{{\partial \mathcal{J}}}{{\partial \tau }}{\mathbf{q}} + \mathcal{J}\frac{{\partial {\mathbf{q}}}}{{\partial \tau }},\boldsymbol \phi } \right) + \frac{1}{2}\sum\limits_{i = 1}^3 {\left( {\frac{\partial\left( {{\mathcal{A}^i}{\mathbf{q}}} \right) }{{\partial {\xi ^i}}},\boldsymbol \phi } \right)}  + \frac{1}{2}\left\{ {\sum\limits_{i = 1}^3 {\left( {\frac{\partial{{\mathcal{A}^i}}}{{\partial {\xi ^i}}}{\mathbf{q}},\boldsymbol \phi } \right)}  + \sum\limits_{i = 1}^3 {\left( {{\mathcal{A}^i}\frac{{\partial {\mathbf{q}}}}{{\partial {\xi ^i}}},\boldsymbol \phi } \right)} } \right\}  = 0.
 \end{equation}
Next we integrate terms that have derivatives of the solution by parts. Note that the coefficient matrices are symmetric so that we can write
\begin{equation}
\begin{split}
\frac{1}{2}\left( {\frac{{\partial \mathcal{J}{\mathbf{q}}}}{{\partial \tau }} + \frac{{\partial \mathcal{J}}}{{\partial \tau }}{\mathbf{q}} + \mathcal{J}\frac{{\partial {\mathbf{q}}}}{{\partial \tau }},\boldsymbol\phi } \right) + {\left. {{{\left( {\vec {\mathcal{A}} \cdot \hat n{\mathbf{q}}} \right)}^T}\boldsymbol\phi } \right|_{\partial E }} &- \frac{1}{2}\left\{ {\sum\limits_{i = 1}^3 {\left( {{\mathcal{A}^i}{\mathbf{q}},\frac{\partial \boldsymbol\phi}{{\partial {\xi ^i}}} } \right)} } \right\} \\& + \frac{1}{2}\left\{ {\sum\limits_{i = 1}^3 {\left( {\frac{\partial {{\mathcal{A}^i}}  }{{\partial {\xi ^i}}}{\mathbf{q}},\boldsymbol\phi } \right)}  - \sum\limits_{i = 1}^3 {\left( {{\mathbf{q}},\frac{\partial \left( {{\mathcal{A}^i}\boldsymbol\phi } \right)}{{\partial {\xi ^i}}}} \right)} } \right\} = 0,
 \label{eq:SkewWeakForm}
 \end{split}
\end{equation}
where we have introduced the shorthand notation
\begin{equation}
{\left. {{{\left( {\vec{\mathcal{A}}  \cdot \hat n{\mathbf{q}}} \right)}^T}\boldsymbol \phi } \right|_{\partial E}} \equiv \sum\limits_{s = 1}^6 {\left\{ {\int_{fac{e^s}} {{{\left( {\vec{\mathcal{A}}  \cdot {{\hat n}^s}{\mathbf{q}}} \right)}^T}\boldsymbol \phi \,d{S^s}} } \right\}},
\end{equation}
for the boundary face contributions. Note that ${\vec{\mathcal{A}}  \cdot \hat n^{s}{\mathbf{q}}} = \vec{\mathcal{A}}{\mathbf{q}}\cdot\hat n^{s} = \tilde {\mathbf{f}}\cdot \hat n^{s}$ is the normal flux at the face $s$.
Eq. (\ref{eq:SkewWeakForm}) is the weak form from which we will create our skew-symmetric approximation.

To get the approximation, (c.f. \cite{Kopriva:2009nx}) we replace ${\mathbf{q}}$ and $\mathcal{J}$ by polynomial interpolants, the normal boundary and interface fluxes by the normal Riemann (numerical) flux, quadratic quantities by their polynomial interpolant and integrals by Legendre-Gauss-Lobatto quadrature. 

We start by defining the polynomial interpolation operator. Let $\mathbb{P}^{N}$ be the space of polynomials of degree less than or equal to $N$. For some function $v\left(\vec\xi\right)$ defined on the reference element, the interpolant of $v$ through the tensor product of the Legendre-Gauss-Lobatto nodes is
\begin{equation}{\mathbb{I}^N}v\left( {\vec \xi } \right) = \sum\limits_{j,k,l = 0}^N {{v_{jkl}}{\ell _j}\left( \xi  \right){\ell _k}\left( \eta  \right){\ell _l}\left( \zeta  \right)} \in \mathbb{P}^{N},\end{equation}
where the $\ell_{j}$ is the Lagrange interpolating polynomial with nodes at the Legendre-Gauss-Lobatto points and ${v_{jkl}}$ is the value of $v$ at the tensor product of those points. We then approximate 
\begin{equation}\begin{gathered}
 {\mathbf{q}}  \approx {\mathbf{Q}} = \sum\limits_{j,k,l = 0}^N {{{\mathbf{Q}}_{jkl}}{\ell _j}\left( \xi  \right){\ell _k}\left( \eta  \right){\ell _l}\left( \zeta  \right)}  \in {\mathbb{P}^N}, \hfill \\
  \mathcal{J} \approx J = \sum\limits_{j,k,l = 0}^N {{J_{jkl}}{\ell _j}\left( \xi  \right){\ell _k}\left( \eta  \right){\ell _l}\left( \zeta  \right)}  \in {\mathbb{P}^N}, \hfill \\
 {{ {\tilde{\mathbf{f}}}}^i}  \approx  {{ \tilde{\mathbf{F}}}^i}= {\mathbb{I}^N}\left( {{\mathbb{I}^N}\left( {{\mathcal{A}^i}} \right){\mathbf{Q}}} \right) = \sum\limits_{j,k,l = 0}^N {\mathcal{A}_{jkl}^i{{{\mathbf{Q}}}_{jkl}}{\ell _j}\left( \xi  \right){\ell _k}\left( \eta  \right){\ell _l}\left( \zeta  \right)}  \in {\mathbb{P}^N}. \hfill \\ 
\end{gathered} \end{equation} We also define the discrete inner product of two polynomials, $\mathbf{ U},\mathbf{ V}\in\mathbb{P}^{N}$
as
\begin{equation}
{\left( {\mathbf U,\mathbf V} \right)_N} = \sum\limits_{j,k,l = 0}^N {\mathbf{ U}_{jkl}^T{{\mathbf{ V}}_{jkl}}{W_j}{W_k}{W_l}}  \equiv \sum\limits_{j,k,l = 0}^N {\mathbf{ U}_{jkl}^T{{\mathbf{ V}}_{jkl}}{W_{jkl}}},
 \end{equation}
where the singly subscripted $W$'s are the one-dimensional Legendre-Gauss-Lobatto quadrature weights and the triply subscripted is the product of the three. We use a similar notation for integrals, where we add a subscript $N$ to denote quadrature. 

We note two facts \cite{CHQZ:2006} about the Legendre-Gauss-Lobatto quadrature
that we will use later. The first is the exactness of the quadrature,
\begin{equation}{\left( {\mathbf{U},\mathbf{V}} \right)_N} = \left( {\mathbf{U},\mathbf{V}} \right)\quad \forall \;{{\mathbf{U}}^T}\mathbf{ V} \in {\mathbb{P}^{2N - 1}}.\end{equation}
The second is that for some function $\mathbf{g}\left(\vec \xi\right)$,
\begin{equation}
{\left( {{\mathbb{I}^N}(\mathbf{ g}),\mathbf{ V}} \right)_N} = {\left( {\mathbf{g},\mathbf{ V}} \right)_N}\quad \forall \;\vec V \in {\mathbb{P}^N},
\label{eq:InterpIP}
\end{equation} which can be seen directly from the definition.

For the numerical flux (Riemann solver) we use the upwind ($\lambda = 1$) or central ($\lambda = 0$)
\begin{equation}
{ {\mathbf{F}}^*}\left( {{{{\mathbf{Q}}}^L},{{{\mathbf{Q}}}^R}};\hat n \right) = \frac{1}{2}\left\{ {\tilde {\mathbf{F}}\left( {{{{\mathbf{Q}}}^L}} \right) \cdot \hat n + \tilde {\mathbf{F}}\left( {{{{\mathbf{Q}}}^R}} \right) \cdot \hat n} \right\} - \lambda \frac{{\left| \mathbb{I}^{N}{\left(\vec {\mathcal{A}}\right) \cdot \hat n} \right|}}{2}\left\{ {{{{\mathbf{Q}}}^R} - {{{\mathbf{Q}}}^L}} \right\},
\label{eq:RiemannSolver}
\end{equation}
which provides a unique flux given a left ${{{\mathbf{Q}}}^L}$ and right ${{{\mathbf{Q}}}^R}$ state relative to the (normal) vector $\hat n$.

We then make the substitutions of the approximations into the continuous weak form to get the discrete weak form
\begin{equation}
\begin{split}
\frac{1}{2}{\left( {\frac{{\partial {\mathbb{I}^N}\left( {J{\mathbf{Q}}} \right)}}{{\partial \tau }} + \frac{{\partial J}}{{\partial \tau }}{\mathbf{Q}} + \mathbb{I}^{N}\left\{J\frac{{\partial {\mathbf{Q}}}}{{\partial \tau }}\right\},\boldsymbol\phi } \right)_N} &+ {\left. {{{\left( {{{ {\mathbf{F}}}^*}} \right)}^T}\boldsymbol\phi } \right|_{\partial E,N }} - \frac{1}{2}\left\{ {\sum\limits_{i = 1}^3 {{{\left( {\tilde {\mathbf{F}}^{i},\frac{\partial \boldsymbol\phi}{{\partial {\xi ^i}}} } \right)}_N}} } \right\} \\& + \frac{1}{2}\left\{ {\sum\limits_{i = 1}^3 {{{\left( {\frac{{\partial {\mathbb{I}^N\left(\mathcal{A}^i\right)}}}{{\partial {\xi ^i}}}{\mathbf{Q}},\boldsymbol\phi } \right)}_N}}  - \sum\limits_{i = 1}^3 {\left( {{\mathbf{Q}},\frac{{\partial {\mathbb{I}^N}\left( {{\mathbb{I}^N\left(\mathcal{A}^i\right)}\boldsymbol\phi } \right)}}{{\partial {\xi ^i}}}} \right)_{N}} } \right\} = 0.
\end{split}\end{equation}
Finally, we separate out the parts that contribute to the GCL to get a formal statement of the approximate weak form
\begin{equation}
\begin{split}
\frac{1}{2}{\left( {\frac{{\partial {\mathbb{I}^N}\left( {J{\mathbf{q}}} \right)}}{{\partial \tau }} + J\frac{{\partial {\mathbf{Q}}}}{{\partial \tau }},\boldsymbol\phi } \right)_N}{\text{ }} &+ {\left. {{{\left( {{{ {\mathbf{F}} }^*}} \right)}^T}\boldsymbol\phi } \right|_{\partial E,N}} - \frac{1}{2}\sum\limits_{i = 1}^3 {{{\left( {{{\tilde {\mathbf{F}} }^i},\frac{{\partial \boldsymbol\phi }}{{\partial {\xi ^i}}}} \right)}_N}}  \\&+ \frac{1}{2}{\sum\limits_{i = 1}^3 {\left( {\frac{{\partial {\mathbb{I}^N}\left( {{{\tilde A}^i}} \right)}}{{\partial {\xi ^i}}}{\mathbf{Q}},\boldsymbol\phi } \right)} _N} - \frac{1}{2}\sum\limits_{i = 1}^3 {{{\left( {{\mathbf{Q}},\frac{{\partial {\mathbb{I}^N}\left( {{\mathbb{I}^N}\left( {{\mathcal{A}^i}} \right)\boldsymbol\phi } \right)}}{{\partial {\xi ^i}}}} \right)}_N}}  \\&+ \frac{1}{2}{\left( {\left[ {\frac{{\partial J}}{{\partial \tau }} - \sum\limits_{i = 1}^3 {\frac{{\partial {\mathbb{I}^N}\left( {\mathcal{J}{{\vec a}^i} \cdot {{\vec x}_\tau }} \right)}}{{\partial {\xi ^i}}}} } \right]{\mathbf{Q}},\boldsymbol\phi } \right)_N} = 0.
\end{split}
\label{eq:SSDGSEMWithGCL}
\end{equation}

\noindent{\bf Reduction to a Static Domain}\\\\
Before moving on, we note that the approximation (\ref{eq:SSDGSEMWithGCL}) is identical to the approximation in \cite{kopriva2015} when the elements do not move. In the static case, time derivatives of the mesh and Jacobian vanish. Next, if we expand the quadratures,
\begin{equation}\frac{1}{2}{\left( {\frac{{\partial {\mathbb{I}^N}\left( {J{\mathbf{Q}}} \right)}}{{\partial \tau }} + J\frac{{\partial {\mathbf{Q}}}}{{\partial \tau }},\boldsymbol\phi } \right)_N} = \frac{1}{2}\sum\limits_{j,k,l=0}^N {{{W}_{jkl}}\boldsymbol\phi \left( {{{\vec \xi }_{jkl}}} \right)\left\{ {\frac{{\partial {J_{jkl}}{{{\mathbf{Q}}}_{jkl}}}}{{\partial \tau }} + {J_{jkl}}\frac{{\partial {{{\mathbf{Q}}}_{jkl}}}}{{\partial \tau }}} \right\}}  = {\left( {J\frac{{\partial {\mathbf{Q}}}}{{\partial \tau }},\boldsymbol\phi } \right)_N}.
\label{eq:timeDerivativeEquivalence}
\end{equation}
Finally, since $\vec x_{\tau} = 0$, $\mathcal{A}^{i} = \tilde A^{i}$. Then for a static domain,
\begin{equation}
\begin{split}
{\left( {J\frac{{\partial {\mathbf{Q}}}}{{\partial \tau }},\boldsymbol \phi } \right)_N} + {\left. {{{\left( {{{ {\mathbf{F}} }^*}} \right)}^T}\boldsymbol \phi } \right|_{\partial E,N}} &- \frac{1}{2}\sum\limits_{i = 1}^3 {{{\left( {{{ \tilde{\mathbf{F}} }^i},\frac{\partial \boldsymbol \phi }{{\partial {\xi ^i}}}} \right)}_N}}  \\&+ \frac{1}{2}\left\{ {{{\left( {\sum\limits_{i = 1}^3 {\frac{{\partial {\mathbb{I}^N}\left( {{{\tilde A}^i}} \right)}}{{\partial {\xi ^i}}}} {\mathbf{Q}},\boldsymbol \phi } \right)}_N} - \sum\limits_{i = 1}^3 {{{\left( {{\mathbf{Q}},\frac{{\partial {\mathbb{I}^N}\left( {{{\tilde A}^i}\boldsymbol \phi } \right)}}{{\partial {\xi ^i}}}} \right)}_N}} } \right\} = 0.
\end{split}
\end{equation}

\subsection{Approximation of the GCL}\label{sec:GCLApprox}

The Geometric Conservation Law (\ref{eq:GCL}) can be written as
\begin{equation}
	\dot{ \mathcal{J}} + \nabla_{\xi}  \cdot \vec \psi  = 0,
\end{equation}
where
\begin{equation}\vec \psi  = -\sum\limits_{i=1}^3 {\mathcal{J}{{\vec a}^i} \cdot {{\vec x}_\tau }{{\hat \xi }^i}}. \end{equation}
For convenience, we approximate this with a DGSEM approximation simultaneously with the solution, so we write a weak form of the GCL as
\begin{equation}
	\left( {{\mathcal{J}_\tau },\phi } \right) + \left( {\nabla_{\xi} \cdot\vec \psi ,\phi } \right) = 0.
\end{equation}
Integrating by parts (to put it into the same equation form as the solution),
\begin{equation}
	\left( {{\mathcal{J}_\tau },\phi } \right) + {\left. {\vec \psi  \cdot {{\hat n} }\phi } \right|_{\partial E} } - \left( {\vec \psi ,\nabla_{\xi} \phi } \right) = 0.
\label{eq:GCLToApproximate}
\end{equation}

We now approximate (\ref{eq:GCLToApproximate}). This means that we 
approximate $\mathcal{J}\approx J\in \mathbb{P}^{N}$ and $\vec \Psi  =  - {\mathbb{I}^N}\left( {\mathcal{J}{{\vec a}^i} \cdot {{\vec x}_\tau }} \right) \in {\mathbb{P}^N}$. Note that by (\ref{eq:contravaraiantBasis}), the flux function, $\vec\Psi$, is actually independent of the Jacobian and is dependent only on the current mesh and its velocity. As a result, (\ref{eq:GCLToApproximate}) doesn't describe a PDE for the Jacobian, but rather is an ODE. Following the recipe above, we replace inner products by Legendre-Gauss-Lobatto quadrature. Formally, we would also replace the boundary term by a Riemann solver. However, the normals (for a conforming mesh) and the mesh velocity are continuous at the faces, so we can simply use the computed values there.
The approximation of the Jacobian is therefore 
\begin{equation}
{\left( {\dot J,\phi } \right)_N} + {\left. {\vec \Psi  \cdot {{\hat n} }\phi } \right|_{\partial E,N}} - {\left( {\vec \Psi ,\nabla \phi } \right)_N} = 0.
\label{eq:DGCL1}
\end{equation}

Furthermore, the discrete inner product satisfies the summation-by-parts (SBP) property \cite{gassner2010}. For some $\tilde {\mathbf{F}}$ and $\phi$,
\begin{equation}
	{\left( {\nabla  \cdot \tilde {\mathbf{F}},\phi } \right)_N} = {\left. {\phi \tilde {\mathbf{F}} \cdot {{\hat n} }} \right|_{\partial E,N}} - {\left( {\tilde {\mathbf{F}},\nabla \phi } \right)_N},
\label{eq:SumByParts}
\end{equation}
so with continuity at the boundaries for $\vec\Psi\cdot\hat n$, the approximation (\ref{eq:DGCL1}) is algebraically equivalent to
\begin{equation}{\left( {\dot J,\phi } \right)_N} + {\left( {\nabla  \cdot \vec \Psi ,\phi } \right)_N} = 0.
\label{eq:DGCL2}
\end{equation}
Finally, we can combine the two inner products to get the equivalent statement
\begin{equation}{\left( {\dot J + \nabla  \cdot \vec \Psi ,\phi } \right)_N} = 0.
\label{eq:DGCL3}
\end{equation}

We now show that the last term in (\ref{eq:SSDGSEMWithGCL}), which contains the GCL, vanishes if we compute the Jacobian using the DG approximation. To see this, it is important to note that the weak form is satisfied \emph{pointwise} by using the quadrature. That is, we take the test function to be the tensor product basis, i.e. $\phi=\ell_j(\xi)\ell_{k}(\eta)\ell_{l}(\zeta)$. Using the last form of the approximation, (\ref{eq:DGCL3}), we see that $\dot J$ satisfies
\begin{equation}
	{{\dot J}_{jkl}} + {\left( {\nabla  \cdot \vec \Psi } \right)_{jkl}} = 0,\quad j,k,l = 0,1, \ldots ,N,
\label{eq:PointwiseJDot}
\end{equation}
or
\begin{equation}
{{\dot J}_{jkl}} - {\left( {\sum\limits_{i = 1}^3 {\frac{{\partial {\mathbb{I}^N}\left( {J{{\vec a}^i} \cdot {{\vec x}_\tau }} \right)}}{{\partial {\xi ^i}}}} } \right)_{jkl}} = 0,\quad j,k,l = 0,1, \ldots ,N,
\label{eq:altForm}
\end{equation}
when we expand the definition of $\vec \Psi$.
Eq.  (\ref{eq:altForm}) holds at each Legendre-Gauss-Lobatto node. Therefore, multiplying by the solution vector, quadrature weight and test function at each node, and summing over all nodes,
\begin{equation}
{\left( {\left[ {\frac{{\partial J}}{{\partial \tau }} - \sum\limits_{i = 1}^3 {\frac{{\partial {\mathbb{I}^N}\left( {J{{\vec a}^i} \cdot {{\vec x}_\tau }} \right)}}{{\partial {\xi ^i}}}} } \right]{\mathbf{Q}},\boldsymbol\phi } \right)_N} = 0.
\label{eq:DGGCLEquiv}
\end{equation}
We will call (\ref{eq:DGGCLEquiv}) the \emph{weak discrete geometric conservation law} or WDGCL. It is also equivalent to the other forms above.

\begin{rem} The approximations
 (\ref{eq:DGCL1}) $\Leftrightarrow$ (\ref{eq:DGCL2}) $\Leftrightarrow$ (\ref{eq:DGCL3}) $\Leftrightarrow$ (\ref{eq:DGGCLEquiv}).
The equivalences will be important later when we want derivatives on the test function or not. They say we can use any of the discrete forms of the GCL as is convenient for theory or computation.
\end{rem}

\subsection{The Skew-Symmetric Approximation on Moving Meshes}

Formally and compactly, provided that the discrete metric identities are satisfied \cite{Kopriva:2006er}, the skew-symmetric approximation on moving meshes for the Jacobian and solution is the geometric conservation law
\begin{equation}
	{\left( { J_{\tau},\phi } \right)_N} + {\left. {\vec \Psi  \cdot {{\hat n} }\phi } \right|_{\partial E,N}} - {\left( {\vec \Psi ,\nabla \phi } \right)_N} = 0,
\label{eq:DGGCL}
\end{equation}
and solution approximation
\begin{equation}
\begin{split}
\frac{1}{2}{\left( {\frac{{\partial {\mathbb{I}^N}\left( {J{\mathbf{Q}}} \right)}}{{\partial \tau }} + J\frac{{\partial {\mathbf{Q}}}}{{\partial \tau }},\boldsymbol\phi } \right)_N}{\text{ }} &+ {\left. {{{\left( {{{ {\mathbf{F}} }^*}} \right)}^T}\boldsymbol\phi } \right|_{\partial E,N}} - \frac{1}{2}\sum\limits_{i = 1}^3 {{{\left( {{{ \tilde{\mathbf{F}} }^i},\frac{{\partial \boldsymbol\phi }}{{\partial {\xi ^i}}}} \right)}_N}} \\& - \frac{1}{2}\sum\limits_{i = 1}^3 {{{\left( {{\mathbf{Q}},\frac{{\partial {\mathbb{I}^N}\left( {{\mathbb{I}^N}\left( {{\mathcal{A}^i}} \right)\boldsymbol\phi } \right)}}{{\partial {\xi ^i}}}} \right)}_N}}  + \frac{1}{2}\sum\limits_{i = 1}^3 {{{\left( {\frac{{\partial {\mathbb{I}^N}\left( {{{\tilde A}^i}} \right)}}{{\partial {\xi ^i}}}{\mathbf{Q}},\boldsymbol\phi } \right)}_N}}  = 0,
\end{split}
\label{eq:SSDGSEM}
\end{equation}
where
\begin{itemize}
\item $\tilde {\mathbf{F}} = \mathbb{I}^{N}\left(\mathbb{I}^{N}\left(\vec {\mathcal{A}}\right){\mathbf{Q}}\right)  = \sum\limits_{i = 1}^3 {{\tilde{\mathbf{F}}^i} {{\hat \xi }^i}} $.
\item $\vec {\mathcal{A}} = \sum\limits_{i = 1}^3 {{{\mathcal{A}}^i}{{\hat \xi }^i}}$.
\item ${\mathcal{A}^i} = J{{\vec a}^i} \cdot \left( {\sum\limits_{j=1}^3 {{A_j}{{\hat x}_j}}  - {{\vec x}_\tau }I} \right) = {{\tilde A}^i} - J{{\vec a}^i} \cdot {{\vec x}_\tau }I$.
\end{itemize}

\section{Properties of the Skew-Symmetric Approximation}\label{PropertiesOfStable}

We now show that the skew-symmetric DGSEM-ALE approximations (\ref{eq:DGGCL}) and (\ref{eq:SSDGSEM}) are stable, conservative, and preserve a constant state when the $A_{j}$'s are constant.

\subsection{Stability}\label{sec:ProvableStability}

The key feature of the skew-symmetric approximation is that it is stable, which follows if the WDGCL is satisfied as described in Sec. \ref{sec:GCLApprox}.

We first derive a bound on the contribution to the energy of a single element. When we set $\boldsymbol\phi = {\mathbf{Q}}$, the time derivative term in (\ref{eq:SSDGSEM}) is
\begin{equation}
\begin{split}
\frac{1}{2}\frac{d}{{d\tau }}\left\| {{\mathbf{Q}}} \right\|_{J,N}^2 = \frac{1}{2}\frac{d}{{d\tau }}{\left( {J{\mathbf{Q}},{\mathbf{Q}}} \right)_N} &= \frac{1}{2}\frac{d}{{d\tau }}\sum\limits_{j,k,l=0}^N {W_{jkl}J_{jkl}{\mathbf{Q}}_{jkl}^T{{{\mathbf{Q}}}_{jkl}}}  \\&= \frac{1}{2}\sum\limits_{j,k,l=0}^N {W_{jkl}\left\{ {\frac{d}{{d\tau }}\left( {{J_{jkl}}{\mathbf{Q}}_{jkl}^T} \right){{\mathbf{Q}}_{jkl}} + {J_{jkl}}\frac{{d{\mathbf{Q}}_{jkl}^T}}{{d\tau }}{{{\mathbf{Q}}}_{jkl}}} \right\}}  \\&= \frac{1}{2}{\left( {\frac{{\partial {\mathbb{I}^N}\left( {J{\mathbf{Q}}} \right)}}{{\partial \tau }} + J\frac{{\partial {\mathbf{Q}}}}{{\partial \tau }},{\mathbf{Q}} } \right)_N}
\end{split}.
\end{equation}
Therefore,
\begin{equation}
\frac{1}{2}\frac{d}{{d\tau }}\left\| {{\mathbf{Q}}} \right\|_{J,N}^2 + {\left. {{{\left( {{{ {\mathbf{F}}}^*}} \right)}^T}{\mathbf{Q}} } \right|_{\partial E,N }} - \frac{1}{2}\sum\limits_{i = 1}^3 {{{\left( { \tilde{\mathbf{F}}^{i},\frac{\partial }{{\partial {\xi ^i}}}{\mathbf{Q}}} \right)}_N}}  - \frac{1}{2}\sum\limits_{i = 1}^3 {\left( {{\mathbf{Q}},\frac{{\partial {\mathbb{I}^N}\left( {{\mathbb{I}^N\left(\mathcal{A}^i\right)}{\mathbf{Q}}} \right)}}{{\partial {\xi ^i}}}} \right)_{N}}+ \frac{1}{2}\sum\limits_{i = 1}^3 {{{\left( {\frac{{\partial {\mathbb{I}^N}\left( {{{\tilde A}^i}} \right)}}{{\partial {\xi ^i}}}{\mathbf{Q}},{\mathbf{Q}} } \right)}_N}} = 0.
\label{eq:Stab1}
\end{equation}
We then apply summation-by-parts (\ref{eq:SumByParts}) to the first sum in (\ref{eq:Stab1}) to move the derivative onto the flux
\begin{equation}
\begin{split}
\frac{1}{2}\frac{d}{{d\tau }}\left\| {{\mathbf{Q}}} \right\|_{J,N}^2 + {\left. {{{\left( {{{ {\mathbf{F}}}^*} - \frac{1}{2}\tilde {\mathbf{F}}\cdot\hat n} \right)}^T}{\mathbf{Q}} } \right|_{\partial E,N }} + \frac{1}{2}\sum\limits_{i = 1}^3 {{{\left( {\frac{\partial }{{\partial {\xi ^i}}} \tilde{\mathbf{F}}^{i},{\mathbf{Q}}} \right)}_N}}  &- \frac{1}{2}\sum\limits_{i = 1}^3 {\left( {{\mathbf{Q}},\frac{{\partial {\mathbb{I}^N}\left( {{\mathbb{I}^N\left(\mathcal{A}^i\right)}{\mathbf{Q}}} \right)}}{{\partial {\xi ^i}}}} \right)_{N}}\\&+ \frac{1}{2}\sum\limits_{i = 1}^3 {{{\left( {\frac{{\partial {\mathbb{I}^N}\left( {{{\tilde A}^i}} \right)}}{{\partial {\xi ^i}}}{\mathbf{Q}},{\mathbf{Q}} } \right)}_N}}  = 0.
\end{split}
\end{equation}
Since $\mathbb{I}^{N}\left(\mathbb{I}^{N}\left(\mathcal{A}^{i}\right){\mathbf{Q}}\right) = \tilde {\mathbf{F}}^{i}$, the two volume terms cancel, leaving
\begin{equation}\frac{1}{2}\frac{d}{{d\tau }}\left\| {{\mathbf{Q}}} \right\|_{J,N}^2 + {\left. {{{\left( {{{ {\mathbf{F}}}^*} - \frac{1}{2}\tilde {\mathbf{F}} \cdot {{\hat n} }} \right)}^T}{\mathbf{Q}}} \right|_{\partial E,N}} = -\frac{1}{2}\sum\limits_{i = 1}^3 {{{\left( {\frac{{\partial {\mathbb{I}^N}\left( {{{\tilde A}^i}} \right)}}{{\partial {\xi ^i}}}{\mathbf{Q}},{\mathbf{Q}} } \right)}_N}} .\end{equation}

If the contravariant coefficient matrices are sufficiently smooth so that the derivatives of the interpolants can be bounded, and if the interpolant of the Jacobian is bounded away from zero \cite{Gassner:2013uq}, then
\begin{equation}
 - \sum\limits_{i = 1}^3 {{{\left( {\frac{{\partial {\mathbb{I}^N}\left( {{{\tilde A}^i}} \right)}}{{\partial {\xi ^i}}}{\mathbf{Q}},{\mathbf{Q}}} \right)}_N}}  \leqslant 2\hat\gamma \left\| {\mathbf{Q}} \right\|_{J,N}^2,
\end{equation}
where
\begin{equation}
\hat \gamma  = \frac{1}{2}\mathop {\max }\limits_{\vec \xi  \in E} \left| \frac{1}{J}{\frac{{\partial {\mathbb{I}^N}\left( {{{\tilde A}^i}} \right)}}{{\partial {\xi ^i}}}} \right|.
\end{equation}

The total energy change is found by summing over all elements. When summed over all elements, the boundary terms along internal faces combine, whereas the boundary terms along physical boundaries do not. Let ${\mathbf{Q}}^{m}$ be the solution on the $m$-th element, $e^{m}$. If we call the interior face contributions $\Sigma_{I}$ and the boundary contributions $\Sigma_{B}$, then
\begin{equation}
\frac{1}{2}\frac{d}{{d\tau }}\sum\limits_{m = 1}^{{N_{el}}} {\left\| {{{\mathbf{Q}}^m}} \right\|_{J,N}^2}  + {\Sigma _I} + {\Sigma _B} \le \frac{{2\hat\gamma }}{2}\sum\limits_{m = 1}^{{N_{el}}} {\left\| {{{\mathbf{Q}}^m}} \right\|_{J,N}^2} .\end{equation}

We now compute the boundary contributions. The external boundary contributions are
\begin{equation}{\Sigma _B} = \sum\limits_{boundary\;faces} {\left\{ {\sum\limits_{r,s = 0}^N {{W_r}{W_s}\left( {{{\mathbf{F}}^*} - \frac{1}{2}\tilde {\mathbf{F}} \cdot \hat n} \right)_{\vec m(r,s)}^T{{{\mathbf{Q}}}_{\vec m(r,s)}}} } \right\}} ,
\label{eq:sigmaB}
\end{equation}
where we use the subscript ``${\vec m(r,s)}$'' to represent the appropriate nodal value on that face. For instance, if the face is the right face of the reference hexahedron then ${\mathbf{Q}}_{\vec m(r,s)}={\mathbf{Q}}_{N,r,s}$.
At each nodal face point, the left state is the computed solution value, and the right state is taken from the exterior of the domain, i.e., from the boundary condition, which we denote as in Sec. \ref{sec:WellPosedness} by ${\mathbf{Q}}_{\infty}$. At each nodal point along a boundary surface,
\begin{equation}
\begin{split}
{\mathbf{F}}^* - \frac{1}{2}\tilde {\mathbf{F}} \cdot \hat n &= \frac{1}{2}\left\{ {\tilde {\mathbf{F}}\left( {{\mathbf{Q}}} \right) \cdot \hat n + \tilde {\mathbf{F}}\left( {{{{\mathbf{Q}}}_\infty }} \right) \cdot \hat n} \right\} - \lambda \frac{{\left| {\mathbb{I}^{N}\left(\vec{\mathcal{A}}\right)  \cdot \hat n} \right|}}{2}\left\{ {{{{\mathbf{Q}}}_\infty } - {\mathbf{Q}}} \right\} - \frac{1}{2}\tilde {\mathbf{F}}\left( {{\mathbf{Q}}} \right) \cdot \hat n \\&= \frac{1}{2}\tilde {\mathbf{F}}\left( {{{{\mathbf{Q}}}_\infty }} \right) \cdot \hat n - \lambda \frac{{\left| {\mathbb{I}^{N}\left(\vec{\mathcal{A}}\right)  \cdot \hat n} \right|}}{2}\left\{ {{{{\mathbf{Q}}}_\infty } - {\mathbf{Q}}} \right\}.
\end{split}
\end{equation}
To guarantee the right kind of bound and therefore stability, we use the upwind solver ($\lambda=1$) at the physical boundaries. For convenience, let us define the intermediate matrix value $A \equiv {\mathbb{I}^N}\left( \vec{\mathcal{A}} \right) \cdot \hat n$. Then the contribution from each boundary point in 
(\ref{eq:sigmaB}) is
\begin{equation}{\left( {{\mathbf{F}}^* - \frac{1}{2}\tilde {\mathbf{F}} \cdot \hat n} \right)^T}{\mathbf{Q}} = \frac{1}{2}{\left\{ {\left( {A - \left| A \right|} \right){{{\mathbf{Q}}}_\infty } + \left| A \right|{\mathbf{Q}}} \right\}^T}{\mathbf{Q}} = \frac{1}{2}\left\{ {{{{\mathbf{Q}}}^T}\left| A \right|{\mathbf{Q}} + 2{{{\mathbf{Q}}}^{T}_\infty }{A^ - }{\mathbf{Q}}} \right\} = \frac{1}{2}\left\{ {{{{\mathbf{Q}}}^T}\left| A \right|{\mathbf{Q}} - 2{{{\mathbf{Q}}}^{T}_\infty }\left| {{A^ - }} \right|{\mathbf{Q}}} \right\},\end{equation}
where $A^{-} = (A -|A|)/2<0$. Since $|A| = A^{+} - A^{-}$, 
\begin{equation}{\left( {{{ {\mathbf{F}} }^*} - \frac{1}{2} \tilde{\mathbf{F}}  \cdot \hat n} \right)^T}{\mathbf{Q}} = \frac{1}{2}\left\{ {{{\mathbf{Q}}^T}{A^ + }{\mathbf{Q}} - {{\mathbf{Q}}^T}{A^ - }{\mathbf{Q}} - 2{\mathbf{Q}}_\infty ^T\left| {{A^ - }} \right|{\mathbf{Q}}} \right\} = \frac{1}{2}{{\mathbf{Q}}^T}{A^ + }{\mathbf{Q}} + \frac{1}{2}\left\{ {{{\mathbf{Q}}^T}\left| {{A^ - }} \right|{\mathbf{Q}} - 2{\mathbf{Q}}_\infty ^T\left| {{A^ - }} \right|{\mathbf{Q}}} \right\}.\end{equation}
We then complete the square on the term in braces to get (c.f. \cite{Gassner:2013uq})
\begin{equation}{\left( {{{ {\mathbf{F}} }^*} - \frac{1}{2} \tilde{\mathbf{F}}  \cdot \hat n} \right)^T}{\mathbf{Q}} = \frac{1}{2}{{\mathbf{Q}}^T}{A^ + }{\mathbf{Q}} + \frac{1}{2}\left\| {\sqrt {\left| {{A^ - }} \right|} {\mathbf{Q}} - \sqrt {\left| {{A^ - }} \right|} {{\mathbf{Q}}_\infty }} \right\|_2^2 - \frac{1}{2}{\mathbf{Q}}_\infty ^T\left| {{A^ - }} \right|{{\mathbf{Q}}_\infty },\end{equation}
where 
${\left\|  \cdot  \right\|_2}$ is the Euclidean 2-norm. Then the boundary face contributions are
\begin{equation}{\Sigma _B} = \frac{1}{2}\sum\limits_{boundary\;faces} {\left\{ {\sum\limits_{r,s = 0}^N {{W_r}{W_s}\left( {\frac{1}{2}{{\mathbf{Q}}^T}{A^ + }{\mathbf{Q}} + \frac{1}{2}\left\| {\sqrt {\left| {{A^ - }} \right|} {\mathbf{Q}} - \sqrt {\left| {{A^ - }} \right|} {{\mathbf{Q}}_\infty }} \right\|_2^2 - \frac{1}{2}{\mathbf{Q}}_\infty ^T\left| {{A^ - }} \right|{{\mathbf{Q}}_\infty }} \right)_{\vec m(r,s)}} } \right\}}. \end{equation}
\begin{rem}
Since the matrix $A$ is symmetric,
\begin{equation}\left\| {\sqrt {\left| {{A^ - }} \right|} {\mathbf{Q}} - \sqrt {\left| {{A^ - }} \right|} {{\mathbf{Q}}_\infty }} \right\|_2^2 = {\left( {{\mathbf{Q}} - {{\mathbf{Q}}_\infty }} \right)^T}\left| {{A^ - }} \right|\left( {{\mathbf{Q}} - {{\mathbf{Q}}_\infty }} \right) \geqslant 0.\end{equation}
\end{rem}
The interior boundary terms have contributions from both the left and the right of a face. When we account for the normals pointing in opposite directions, the face terms are
\begin{equation}{\Sigma _I} = \frac{1}{2}\sum\limits_{interior\;faces} {\left\{ {\sum\limits_{r,s = 0}^N {{W_r}{W_s}\left( {F_{\vec m(r,s)}^{*T}\jump{{{{{\mathbf{Q}}}_{\vec m(r,s)}}}} - \frac{1}{2}\jump{ {\tilde {\mathbf{F}}_{\vec m(r,s)}^T \cdot \hat n{{{\mathbf{Q}}}_{\vec m(r,s)}}} }} \right)} } \right\}}, \end{equation}
where $\jump{\cdot}$ represents the jump in the quantity across the interface. At each face point \cite{Gassner:2013uq}, 
\begin{equation}{\mathbf{F}^{*T}}\jump{ {{\mathbf{Q}}} } - \frac{1}{2}\jump{ {{\tilde{\mathbf{F}}^T}{\mathbf{Q}}} } = \frac{\lambda }{2}{\jump{ {{\mathbf{Q}}} }^T}\left| A \right|\jump{ {{\mathbf{Q}}} } \geqslant 0.\end{equation}
Therefore, $\Sigma_{I}\ge 0$ and
\begin{equation}
\begin{split}
\frac{1}{2}\frac{d}{{d\tau }}\sum\limits_{m = 1}^{{N_{el}}} {\left\| {{{\mathbf{Q}}^m}} \right\|_{J,N}^2}  \leqslant  &- \frac{1}{2}\sum\limits_{boundary\;faces} {\left\{ {\sum\limits_{r,s = 0}^N {{W_r}{W_s}\left( {{{\mathbf{Q}}^T}{A^ + }{\mathbf{Q}} + \left\| {\sqrt {{{\left| {{A^ - }} \right|}_{\vec m(r,s)}}} {{\mathbf{Q}}_{\vec m(r,s)}} - \sqrt {{{\left| {{A^ - }} \right|}_{\vec m(r,s)}}} {{\mathbf{Q}}_\infty }} \right\|_2^2} \right)} } \right\}}  \\&+ \frac{1}{2}\sum\limits_{boundary\;faces} {\left\{ {\sum\limits_{r,s = 0}^N {{W_r}{W_s}{\mathbf{Q}}_\infty ^T{{\left| {{A^ - }} \right|}_{\vec m(r,s)}}{{\mathbf{Q}}_\infty }} } \right\}}  + \hat \gamma \sum\limits_{m = 1}^{{N_{el}}} {\left\| {{{\mathbf{Q}}^m}} \right\|_{J,N}^2} .
\end{split}
\label{eq:GlobalEnergyODE}
\end{equation}

Let us now define the norm over the entire domain to be
\begin{equation}{\left\| {{\mathbf{Q}}} \right\|^{2}_{J,N}} \equiv \sum\limits_{m = 1}^{{N_{el}}} {\left\| {{{{\mathbf{Q}}}^m}} \right\|_{J,N}^2}, \end{equation}
and integrate both sides of (\ref{eq:GlobalEnergyODE}) with respect to time. Then the norm of the solution can be bounded in terms of the initial and boundary values by (c.f. (\ref{eq:StronglyWellPosed}))
\begin{equation}
\begin{split}
\left\| {{\mathbf{Q}}(T)} \right\|_{J,N}^2 &+ {e^{2\hat \gamma T}}\int_0^T {\sum\limits_{boundary\;faces} {\left\{ {\sum\limits_{r,s = 0}^N {{W_r}{W_s}\left( {{{\mathbf{Q}}^T}{A^ + }{\mathbf{Q}} + \left\| {\sqrt {{{\left| {{A^ - }} \right|}_{\vec m(r,s)}}}\left( {{\mathbf{Q}}_{\vec m(r,s)}} - {{\mathbf{Q}}_\infty }\right)} \right\|_2^2} \right)} } \right\}} dt} \\& \leqslant {e^{2\hat \gamma T}}\left\{ {\left\| {{\mathbf{Q}}(0)} \right\|_{J,N}^2 + \int_0^T {\sum\limits_{boundary\;faces} {\left\{ {\sum\limits_{r,s = 0}^N {{W_r}{W_s}{\mathbf{Q}}_\infty ^T{{\left| {{A^ - }} \right|}_{\vec m(r,s)}}{{\mathbf{Q}}_\infty }} } \right\}} dt} } \right\}.
\end{split}
 \end{equation}
 Finally, we note that the normed term on the left, being positive, represents additional dissipation along the incoming characteristics due to the upwind Riemann solver, which goes to zero as the solution converges. Since it is positive, we can also bound the solution by
 \begin{equation}
 \begin{split}
 \left\| {{\mathbf{Q}}(T)} \right\|_{J,N}^2 &+ {e^{2\hat \gamma T}}\int_0^T {\sum\limits_{boundary\;faces} {\left\{ {\sum\limits_{r,s = 0}^N {{W_r}{W_s}{{\mathbf{Q}}_{\vec m(r,s)}^T}{A^+_{\vec m(r,s)} }{\mathbf{Q}}_{\vec m(r,s)}} } \right\}} dt}  \\&\leqslant {e^{2\hat \gamma T}}\left\{ {\left\| {{\mathbf{Q}}(0)} \right\|_{J,N}^2 + \int_0^T {\sum\limits_{boundary\;faces} {\left\{ {\sum\limits_{r,s = 0}^N {{W_r}{W_s}{\mathbf{Q}}_\infty ^T{{\left| {{A^ - }} \right|}_{\vec m(r,s)}}{{\mathbf{Q}}_\infty }} } \right\}} dt} } \right\}.
\label{eq:stronglyStable}
 \end{split}\end{equation}
 \begin{rem}
 The bound (\ref{eq:stronglyStable}), which is the discrete equivalent of the strong well-posedness bound (\ref{eq:StronglyWellPosed}), says that the approximation is \emph{strongly stable}.
 \end{rem}
\begin{rem}
Note that the statement (\ref{eq:stronglyStable}) implies that the approximation is also \emph{stable}, for when ${\mathbf{Q}}_{\infty}= 0$,
\begin{equation}
{\left\| {{\mathbf{Q}}(T)} \right\|_{J,N}} \leqslant e^{\hat\gamma T}{\left\| {{\mathbf{Q}}(0)} \right\|_{J,N}}.
\end{equation}
Furthermore, for constant coefficient problems where the discrete metric identities are satisfied, $\hat \gamma =0$ \cite{kopriva2015} so
\begin{equation}
{\left\| {{\mathbf{Q}}(T)} \right\|_{J,N}} \leqslant {\left\| {{\mathbf{Q}}(0)} \right\|_{J,N}}.
\end{equation}
\end{rem}
\begin{rem}
The discrete and continuous norms are equivalent \cite{Gassner:2013uq}. Therefore the continuous norms for the approximate solution satisfy similar bounds.
\end{rem}

\subsection{Conservation}\label{ConservationSection}

To show that the approximation is conservative, we start on a single element and add (\ref{eq:SSDGSEM}) to the WDGCL equivalent form (\ref{eq:DGGCLEquiv}) to get back to (\ref{eq:SSDGSEMWithGCL}). We then combine terms that contribute to the contravariant coefficient matrices and rewrite it here as
\begin{equation}
\begin{split}
\frac{1}{2}{\left( {\frac{{\partial {\mathbb{I}^N}\left( {J{\mathbf{Q}}} \right)}}{{\partial \tau }} + J\frac{{\partial {\mathbf{Q}}}}{{\partial \tau }},\boldsymbol \phi } \right)_N} + &{\left. {{{\left( {{{ {\mathbf{F}}}^*}} \right)}^T}\boldsymbol \phi } \right|_{\partial E,N }} - \frac{1}{2}\sum\limits_{i = 1}^3 {{{\left( {\tilde {\mathbf{F}}^{i},\frac{\partial }{{\partial {\xi ^i}}}\boldsymbol \phi } \right)}_N}}  \\&- \frac{1}{2}\sum\limits_{i = 1}^3 {\left( {{\mathbf{Q}},\frac{{\partial {\mathbb{I}^N}\left( {{\mathcal{A}^i}\phi } \right)}}{{\partial {\xi ^i}}}} \right)_{N}}  + \frac{1}{2}{\left( {\left[ {\frac{{\partial J}}{{\partial \tau }} + \sum\limits_{i = 1}^3 {\frac{{\partial {\mathbb{I}^N}\left( {{\mathcal{A}^i}} \right)}}{{\partial {\xi ^i}}}} } \right]{\mathbf{Q}},\boldsymbol \phi } \right)_N} = 0.
\end{split}
\label{eq:SummedSSDGSEMAndDGCL}\end{equation}
We then let $\boldsymbol\phi = \hat {\mathbf{e}}_{n}$, which is the vector with one in the $n^{th}$ entry and zero otherwise. In other words, we separate out the equations in the system and choose the test function to be one for each in turn.

Gathering the time derivative terms, and with a slight bit of notational abuse to go from the single equation associated with $\hat {\mathbf{e}}_{n}$ to the state vector form that incorporates each $n$,
\begin{equation}\begin{split} 
\frac{1}{2}{\left( {\frac{{\partial {\mathbb{I}^N}\left( {J{\mathbf{Q}}} \right)}}{{\partial \tau }} + J\frac{{\partial {\mathbf{Q}}}}{{\partial \tau }} + \frac{{\partial J}}{{\partial \tau }}{\mathbf{Q}},{{\hat {\mathbf{e}}}_n}} \right)_N} &\rightarrow\frac{1}{2}\sum\limits_{i,j,k = 0}^N {{W_{ijk}}\left\{ {\frac{{\partial {J_{ijk}}{{{\mathbf{Q}}}_{ijk}}}}{{\partial \tau }} + {J_{ijk}}\frac{{\partial {{{\mathbf{Q}}}_{ijk}}}}{{\partial \tau }} + \frac{{\partial {J_{ijk}}}}{{\partial \tau }}{{{\mathbf{Q}}}_{ijk}}} \right\}}  \\&= \sum\limits_{i,j,k = 0}^N {{W_{ijk}}\left\{ {\frac{{\partial {J_{ijk}}{{{\mathbf{Q}}}_{ijk}}}}{{\partial \tau }}} \right\}}  = \frac{d}{{d\tau }}\sum\limits_{i,j,k = 0}^N {{W_{ijk}}{J_{ijk}}{{{\mathbf{Q}}}_{ijk}}} \\& \equiv \frac{d}{d\tau}\int\limits_{E ,N} {J{\mathbf{Q}}\,d\vec \xi },
 \end{split}\end{equation}
 which gives the time rate of change of the total amount of ${\mathbf{Q}}$ in an element.
 
 Next, it is immediate that the term 
\begin{equation}\sum\limits_{i = 1}^3 {{{\left( { \tilde{\mathbf{F}}^{i},\frac{{\partial {{\hat {\mathbf{e}}}_n}}}{{\partial {\xi ^i}}}} \right)}_N}}  = 0.\end{equation}
Finally, we see that the internal flux terms cancel. Since $\mathcal{A}^{i}$ is symmetric,
\begin{equation}
\begin{split}
& - \frac{1}{2}\sum\limits_{i = 1}^3 {\left( {{\mathbf{Q}},\frac{{\partial {\mathbb{I}^N}\left( {{\mathcal{A}^i}{{\hat {\mathbf{e}}}_n}} \right)}}{{\partial {\xi ^i}}}} \right)_{N}}  + \frac{1}{2}{\left( {\left[ {\sum\limits_{i = 1}^3 {\frac{{\partial {\mathbb{I}^N}\left( {{\mathcal{A}^i}} \right)}}{{\partial {\xi ^i}}}} } \right]{\mathbf{Q}},{{\hat{\mathbf{e}}}_n}} \right)_N} \\&=  - \frac{1}{2}\sum\limits_{i = 1}^3 {\left( {{\mathbf{Q}},\frac{{\partial {\mathbb{I}^N}\left( {{\mathcal{A}^i}} \right)}}{{\partial {\xi ^i}}}{{\hat {\mathbf{e}}}_n}} \right)_{N}}  + \frac{1}{2}\sum\limits_{i = 1}^3 {{{\left( {\frac{{\partial {\mathbb{I}^N}\left( {{\mathcal{A}^i}} \right)}}{{\partial {\xi ^i}}}{\mathbf{Q}},{{\hat {\mathbf{e}}}_n}} \right)}_N}}  = 0.
 \end{split}
 \end{equation}
Therefore, we have conservation on each element,
\begin{equation}\frac{d}{{d\tau }}\int\limits_{E,N} {J{\mathbf{Q}}\,d\vec \xi }  =  - { {{{\left( {{{ {\mathbf{F}}}^*}} \right)}}} \bigg|_{\partial E,N}},\end{equation}
which says that the rate of change of the total amount of ${\mathbf{Q}}$ depends only on the net flux on the element faces. Therefore, locally on each element, the solutions satisfy a discrete conservation property and the quantity
\begin{equation}\label{conservedQuantity}
\sum_{i,j,k=0}^N(J\mathbf{Q})_{ijk}W_{ijk},
\end{equation}
is conserved.

The global conservation statement is found when we sum over all the elements. We note that the use of the Riemann solver gives an equal and opposite flux at the faces, so the interior flux contributions cancel exactly. What is left is determined only by the net flux through the boundaries of the domain,
\begin{equation}\frac{d}{{d\tau }}\sum\limits_{m = 1}^{{N_{el}}} {\int\limits_{E,N} {{J^m}{{{\mathbf{Q}}}^m}\,d\vec \xi } }  =  - \sum\limits_{boundary\;faces} {\left\{ {\sum\limits_{r,s = 0}^N {{W_r}{W_s}\tilde {\mathbf{F}}_{\vec m(r,s)}^*} } \right\}} .\end{equation}

\subsection{Free-Stream Preservation}\label{FreeStreamSection}

Finally, we assume that the solution is a constant and show that it remains constant for all time if the coefficient matrices are constant. 

We assume that the metric identities (\ref{eq:MetricIdentities}) hold discretely by construction. That is, we assume the metric terms are computed from the transformation to the reference element so that \cite{Kopriva:2006er}
\begin{equation}
	\sum\limits_{i=1}^3 {\frac{\partial \mathbb{I}^{N}\left({J{{\vec a}^i}}\right)}{{\partial {\xi ^i}}}}  = 0,
\end{equation}
and since the coefficient matrices are constant,
\begin{equation}
	\sum\limits_{i=1}^3 {\frac{{\partial {\mathbb{I}^N}\left( {J{{\vec a}^i} \cdot \sum\limits_{j=1}^3 {{A_j}{{\hat x}_j}} } \right)}}{{\partial {\xi ^i}}}}  = \sum\limits_{i=1}^3 {\frac{{\partial {\mathbb{I}^N}\left( {{{\tilde A}^i}} \right)}}{{\partial {\xi ^i}}}}  = 0.
\label{eq:AtlMetric}
\end{equation}
Eq.  (\ref{eq:altForm}) holds at each node in an element. Therefore, multiplying by the solution vector at each node,
\begin{equation}{{\dot J}_{jkl}}{{{\mathbf{Q}}}_{jkl}} + {\left( {\sum\limits_{i=1}^3 {\frac{{\partial {\Psi ^i}}}{{\partial {\xi ^i}}}} } \right)_{jkl}}{{{\mathbf{Q}}}_{jkl}} + \sum\limits_{i=1}^3 {\frac{{\partial {\mathbb{I}^N}\left( {{{\tilde A}^i}} \right)}}{{\partial {\xi ^i}}}{{{\mathbf{Q}}}_{jkl}}}  = 0,\quad j,k,l= 0,1, \ldots ,N.\end{equation}
We can then combine the sums
\begin{equation}{{\dot J}_{jkl}}{{{\mathbf{Q}}}_{jkl}} + {\left( {\sum\limits_{i=1}^3 {\frac{{\partial {\mathbb{I}^N}\left( {{\mathcal{A}^i}} \right)}}{{\partial {\xi ^i}}}} } \right)_{jkl}}{{{\mathbf{Q}}}_{jkl}} = 0,\quad j,k,l = 0,1, \ldots ,N.\end{equation}
Finally, we multiply by the quadrature weights and the test function at each point and sum over all nodes to see that when the coefficient matrices are constant, and the metric identities are discretely satisfied, then the WDGCL is equivalent to
\begin{equation}\frac{1}{2}{\left( {{J_\tau }{\mathbf{Q}},\boldsymbol \phi } \right)_N} =  - \frac{1}{2}{\left( {\sum\limits_{i = 1}^3 {\frac{{\partial {\mathbb{I}^N}\left( {{\mathcal{A}^i}} \right)}}{{\partial {\xi ^i}}}{\mathbf{Q}}} ,\boldsymbol \phi } \right)_N}.\end{equation} 

We then note that
\begin{equation}\frac{1}{2}{\left( {\frac{{\partial {\mathbb{I}^N}\left( {J{\mathbf{Q}}} \right)}}{{\partial \tau }} + J\frac{{\partial {\mathbf{Q}}}}{{\partial \tau }},\boldsymbol \phi } \right)_N} = {\left( {J\frac{{\partial {\mathbf{Q}}}}{{\partial \tau }},\boldsymbol \phi } \right)_N} + \frac{1}{2}{\left( {{J_\tau }{\mathbf{Q}},\boldsymbol \phi } \right)_N}.\end{equation}
Therefore, (\ref{eq:SSDGSEM}) is equivalent to
\begin{equation}{\left( {J\frac{{\partial {\mathbf{Q}}}}{{\partial \tau }},\boldsymbol \phi } \right)_N}{\text{ }} + {\left. {{{\left( {{{ {\mathbf{F}} }^*}} \right)}^T}\boldsymbol \phi } \right|_{\partial E}} - \frac{1}{2}\sum\limits_{i = 1}^3 {{{\left( {{{\tilde {\mathbf{F}} }^i},\frac{\partial }{{\partial {\xi ^i}}}\boldsymbol \phi } \right)}_N}}  - \frac{1}{2}\sum\limits_{i = 1}^3 {{{\left( {{\mathbf{Q}},\frac{{\partial {\mathbb{I}^N}\left( {{\mathbb{I}^N}\left( {{\mathcal{A}^i}} \right)\boldsymbol \phi } \right)}}{{\partial {\xi ^i}}}} \right)}_N}}  - \frac{1}{2}{\left( {\sum\limits_{i = 1}^3 {\frac{{\partial {\mathbb{I}^N}\left( {{\mathcal{A}^i}} \right)}}{{\partial {\xi ^i}}}{\mathbf{Q}}} ,\boldsymbol \phi } \right)_N} = 0.
\label{eq:SSDGSEMMinusGCL}
\end{equation}

We now set ${\mathbf{Q}} = \mathbf{c}$, where $\mathbf{c}$ is a constant, and show that the time derivative of the solution is zero at each nodal point in an element.
First we use the consistency of the Riemann solver (\ref{eq:RiemannSolver}) to show that it is equivalent to the normal boundary flux. Using the internal and external states as ${\mathbf{Q}}^{L}={\mathbf{Q}}^{R} = \mathbf{c}$,
\begin{equation}
{{ {\mathbf{F}}}^*}\left( {\mathbf{c},\mathbf{c}} \right) = \frac{1}{2}\left\{ {\tilde {\mathbf{F}}\left( \mathbf{c}\right) \cdot \hat n + \tilde {\mathbf{F}}\left( {\mathbf{c}} \right) \cdot \hat n} \right\} = \tilde {\mathbf{F}}\left( {\mathbf{c}} \right) \cdot \hat n.
\end{equation}
Next, we use summation-by-parts on the first sum in (\ref{eq:SSDGSEMMinusGCL}):
\begin{equation}
 - \frac{1}{2}\sum\limits_{i = 1}^3 {{{\left( {\tilde {\mathbf{F}},\frac{\partial }{{\partial {\xi ^i}}}\boldsymbol \phi } \right)}_N}}  =  - \frac{1}{2}{\left. {{{\left( {\tilde {\mathbf{F}} \cdot \hat n} \right)}^T}\boldsymbol \phi } \right|_{\partial E,N}} + \frac{1}{2}\sum\limits_{i = 1}^3 {{{\left( {\frac{\partial }{{\partial {\xi ^i}}}{{ \tilde{\mathbf{F}}}^i},\boldsymbol \phi } \right)}_N}}  =  - \frac{1}{2}{\left. {{{\left( {\tilde {\mathbf{F}} \cdot \hat n} \right)}^T}\boldsymbol \phi } \right|_{\partial E,N}} + \frac{1}{2}\sum\limits_{i = 1}^3 {{{\left( {\frac{{\partial {\mathbb{I}^N}\left( {{\mathcal{A}^i}} \right)}}{{\partial {\xi ^i}}}\mathbf{c},\boldsymbol \phi } \right)}_N}},  \label{eq:SSDGSEMPart}
\end{equation}
and note that the last sum in (\ref{eq:SSDGSEMPart}) will cancel the equivalent term from the GCL in (\ref{eq:SSDGSEMMinusGCL}).

We also use summation-by-parts and (\ref{eq:InterpIP}) on the second sum in (\ref{eq:SSDGSEMMinusGCL}). Since $\mathbf{c}$ is a constant, the volume term vanishes, leaving only the boundary term
\begin{equation}
 - \frac{1}{2}\sum\limits_{i = 1}^3 {{{\left( {\mathbf{ c},\frac{{\partial {\mathbb{I}^N}\left( {{\mathcal{A}^i}\boldsymbol \phi } \right)}}{{\partial {\xi ^i}}}} \right)}_N}}  =  - \frac{1}{2}{\left. {{{\mathbf{ c}}^T}\left( {\vec{ \mathcal{A}} \boldsymbol \phi } \right)} \right|_{\partial E,N}} =  - \frac{1}{2}{\left. {{{\left( {\vec{ \mathcal{A}}  \cdot \hat n\mathbf{ c}} \right)}^T}\boldsymbol \phi } \right|_{\partial E,N}} =  - \frac{1}{2}{\left. {{{\left( {\tilde {\mathbf{F}} \cdot \hat n} \right)}^T}\boldsymbol \phi } \right|_{\partial E,N}}.
 \label{eq:SecondSumResolution}
 \end{equation}

Now we replace the terms in (\ref{eq:SSDGSEMMinusGCL}) with (\ref{eq:SSDGSEMPart}) and (\ref{eq:SecondSumResolution}). Note that the boundary term from the Riemann solver is cancelled by the two terms that come from the summation-by-parts. The volume term from the first sum in (\ref{eq:SSDGSEMMinusGCL}) is cancelled by the sum term in the GCL part. All this leaves is
\begin{equation}
{\left( {J\frac{{\partial {\mathbf{Q}}}}{{\partial \tau }},\boldsymbol \phi } \right)_N}=0.
\end{equation}
To show the final result, we let $\boldsymbol \phi = \ell_{i}\left(\xi\right)\ell_{j}\left(\eta\right)\ell_{k}\left(\zeta\right)\hat{\mathbf{e}}_{n}$, use the discrete orthogonality of the Lagrange interpolating polynomials, and divide by the quadrature weights to see that
\begin{equation}
\frac{{\partial {{{\mathbf{Q}}}_{ijk}}}}{{\partial \tau }} = 0,
\end{equation}
at each nodal point in the mesh, so the approximation preserves a constant state. 

\begin{rem}
The constant solution of the PDE (\ref{eq:OrigConsLaw}) also stays constant if the coefficient matrices are variable but are divergence free, that is,
\begin{equation}
\sum\limits_{i = 1}^3 {\frac{{\partial {A_i}(\vec x)}}{{\partial {x_i}}}}  = 0.
\end{equation}
When the coefficient matrices are variable, it is not immediate that (\ref{eq:AtlMetric}) holds. The interpolation of the product introduces aliasing error when $\vec A$ is not constant, and since the differentiation and interpolation will not commute, there will be a (spectrally) small error in the time derivative of the solution due to the divergence of the metric identities dotted with the vector of coefficient matrices not being zero. 

Several strategies could be used to retain constant state preservation discretely. In some special circumstances it might be possible to write the product as a curl of some quantity. More generally, one can compute the nodal values of the product ${J{{\vec a}^i} \cdot \sum\limits_{j=1}^3 {{A_j}{{\hat x}_j}} }$ by a projection method \cite{Birdsall:1985qv}, which requires the solution of a Poisson problem. For static meshes, the calculation can be done once and stored, but for moving meshes, solving the Poisson problem at each time step will likely be expensive. Another alternative is to use hyperbolic divergence cleaning \cite{Munz:2000fr}. See \cite{Jacobs:2006qd} for a discussion of the projection and hyperbolic divergence cleaning methods applied to discontinuous Galerkin spectral element approximations.
\end{rem}

\section{Implementation}\label{ImplementationSection}

We now re-write the weak form (\ref{eq:SSDGSEM}) into one suitable for computation. For notational convenience, we define
\begin{equation}
\dot {\mathbf{H}} \equiv \frac{1}{2}\left\{ {\frac{{\partial {\mathbb{I}^N}\left( {J{\mathbf{Q}}} \right)}}{{\partial \tau }} + J\frac{{\partial {\mathbf{Q}}}}{{\partial \tau }}} \right\}.
\end{equation}
Then the time derivative of the solution will be computed through
\begin{equation}
{\left( {\dot {\mathbf{H}},\boldsymbol{\phi} } \right)_N} + {\left. {{{\left( {{{ {\mathbf{F}} }^*}} \right)}^T}\boldsymbol{\phi} } \right|_{\partial E,N}} + \frac{1}{2}\sum\limits_{i = 1}^3 {{{\left( {{\tilde{\mathbf{F}}^i},\frac{{\partial \boldsymbol{\phi} }}{{\partial {\xi ^i}}}} \right)}_N}}  - \frac{1}{2}\sum\limits_{i = 1}^3 {{{\left( {{\mathbf{Q}},\frac{{\partial {\mathbb{I}^N}\left( {{\mathbb{I}^N}\left( {{\mathcal{A}^i}} \right)\boldsymbol{\phi} } \right)}}{{\partial {\xi ^i}}}} \right)}_N}}  + \frac{1}{2}\sum\limits_{i = 1}^3 {{{\left( {\frac{{\partial {\mathbb{I}^N}\left( {{{\tilde A}^i}} \right)}}{{\partial {\xi ^i}}}{\mathbf{Q}},\boldsymbol{\phi} } \right)}_N}}=0.
\label{eq:Hdot}
\end{equation}
The GCL computes the time derivative $\dot J$ at each point through
(\ref{eq:DGCL1}). See also \cite{Minoli:2010rt}.

The pointwise definitions of the time derivatives are found by replacing $\boldsymbol\phi$ by the Lagrange basis functions in each equation. In that way, at each node in an element,
\begin{equation}{\left( {\dot {\mathbf{H}},{\ell _j}{\ell _k}{\ell _l}{{\widehat {\mathbf{e}}}_n}} \right)_N} \to {{\dot {\mathbf{H}}}_{jkl}}{W_{jkl}}.\end{equation}

We next find the pointwise representations of the terms on the right. Except for the one half, the first and second terms on the right of (\ref{eq:Hdot}) reduce to the standard DGSEM formula \cite{Kopriva:2009nx}. The volume term is therefore
\begin{equation}\sum\limits_{i = 1}^3 {{{\left( {{{\mathbf{F}}^i},\frac{{\partial\boldsymbol{\phi}}}{{\partial {\xi ^i}}}} \right)}_N}}  =  - {W_{jkl}}\left\{ {\sum\limits_{n = 0}^N {\tilde{\mathbf{F}}_{nkl}^1{{\hat D}_{jn}}}  + \sum\limits_{n = 0}^N {\tilde{\mathbf{F}}_{jnl}^2{{\hat D}_{kn}}}  + \sum\limits_{n = 0}^N {\tilde{\mathbf{F}}_{jkn}^3{{\hat D}_{ln }}} } \right\},
\end{equation}
where
\begin{equation}
	\hat D_{jn}  =  - \dfrac{{D_{nj} W_n }}{{W_j }},
	\label{eq:dhat}
\end{equation}
and	$D_{nj}  = \ell '_j \left( {\xi_n } \right)$, etc. The surface terms are identical to the standard formulation found in \cite{Kopriva:2009nx}. The derivatives of the contravariant coefficient matrices in the last term can be computed by standard matrix-vector multiplication using the standard derivative matrix by
\begin{equation}{\left. {\sum\limits_{i = 1}^3 {\frac{{\partial {\mathbb{I}^N}\left( {{{\tilde A}^i}} \right)}}{{\partial {\xi ^i}}}} } \right|_{mnp}} = \sum\limits_{r = 0}^N {\tilde A_{rnp}^1{D_{mr}}}  + \sum\limits_{r = 0}^N {\tilde A_{mrp}^2{D_{nr}}}  + \sum\limits_{r = 0}^N {\tilde A_{mnr}^3{D_{pr}}}. \end{equation}
Therefore, the last volume term on the right of (\ref{eq:Hdot}) has the form
\begin{equation}
\sum\limits_{i = 1}^3 {{{\left( {\frac{{\partial {\mathbb{I}^N}\left( {{{\tilde A}^i}} \right)}}{{\partial {\xi ^i}}}{\mathbf{Q}},\boldsymbol{\phi} } \right)}_N}}  \to {W_{jkl}}\left\{ {\sum\limits_{r = 0}^N {\tilde A_{rkl}^1{D_{jr}}}  + \sum\limits_{r = 0}^N {\tilde A_{jrl}^2{D_{kr}}}  + \sum\limits_{r = 0}^N {\tilde A_{jkr}^3{D_{lr}}} } \right\}{{\mathbf{Q}}_{jkl}} \equiv {W_{jkl}}{G_{jkl}}{{\mathbf{Q}}_{jkl}}.
\end{equation}
For efficiency, $G_{jkl}$ can be computed once and stored. If the coefficient matrices are constant, and the discrete metric identities are satisfied, it is zero.

The most interesting term essentially calculates the derivatives of the fluxes computed from the test functions rather than the solution values. Because of the tensor product form of the solution, it is enough to derive the representation using a generic one dimensional inner product of the form
\begin{equation}{{{\left( {{{\left( {{\mathbb{I}^N}\left( {A\boldsymbol \phi } \right)} \right)}^\prime },\mathbf U} \right) }_N}}.\end{equation}
First,
\begin{equation}{\mathbb{I}^N}\left( {A\boldsymbol \phi } \right) = \sum\limits_{n = 0}^N {A\left( {{\xi _n}} \right)\boldsymbol \phi \left( {{\xi _n}} \right){\ell _n}\left( \xi  \right)}, \end{equation}
so
\begin{equation}{\left( {{\mathbb{I}^N}\left( {A\boldsymbol \phi } \right)} \right)^\prime } = \sum\limits_{n = 0}^N {A\left( {{\xi _n}} \right)\boldsymbol \phi \left( {{\xi _n}} \right){{\ell '}_n}\left( \xi  \right)}. \end{equation}
For the matrices we consider, $A$ is symmetric, so
\begin{equation}{\left( {{{\left( {{\mathbb{I}^N}\left( {A\boldsymbol \phi } \right)} \right)}^\prime },\mathbf{ U}} \right) _N} = \sum\limits_{m = 0}^N {{W_m}\sum\limits_{n = 0}^N {{{\ell '}_n}\left( \xi _{m} \right){{\left( {A\left( {{\xi _n}} \right)\boldsymbol \phi \left( {{\xi _n}} \right)} \right)}^T}} {{\mathbf{ U}}_m}}  = \sum\limits_{m = 0}^N {{W_m}\sum\limits_{n = 0}^N {{{\ell '}_n}\left( {{\xi _m}} \right){{\boldsymbol \phi }^T}\left( {{\xi _n}} \right)A\left( {{\xi _n}} \right)} {{\mathbf{ U}}_m}}. \end{equation}
Now we take $\boldsymbol \phi = \hat {\mathbf{e}}_{k}\ell_{j}\left(\xi\right)$ to select out each component of the vector,
\begin{equation}{\left( {{{\left( {{I^N}\left( {A\hat {\mathbf{e}}_{k}\ell_{j} } \right)} \right)}^\prime },\mathbf{U}} \right) _N} = \sum\limits_{m = 0}^N {{W_m}\sum\limits_{n = 0}^N {{{\ell '}_n}\left( {{\xi _m}} \right){\ell _j}\left( {{\xi _n}} \right)\hat {\mathbf{e}}_k^TA\left( {{\xi _n}} \right)} {{\mathbf{U}}_m}}. \end{equation}
Using the Kronecker Delta property of the Lagrange interpolating polynomials, only one term of the $n$ sum remains so that
\begin{equation}{\left( {{{\left( {{I^N}\left( {A\hat {\mathbf{e}}_{k}\ell_{j} } \right)} \right)}^\prime },\mathbf{ U}} \right) _N} = \sum\limits_{m = 0}^N {\hat{\mathbf{ e}}_k^TA\left( {{\xi _n}} \right){{\mathbf{ U}}_m}{W_m}{{\ell '}_n}\left( {{\xi _m}} \right)}  = -{W_j}\sum\limits_{m = 0}^N {\hat{\mathbf{ e}}_k^TA\left( {{\xi _j}} \right){{\mathbf{ U}}_m}{{\hat D}_{jm}}}. \end{equation}
Finally, we reconstruct the vector from the components to see that
\begin{equation}{\left( {{{\left( {{I^N}\left( {A\phi } \right)} \right)}^\prime },{\mathbf{U}}} \right)_N} \to  - {W_j}\sum\limits_{m = 0}^N {A\left( {{\xi _j}} \right){{\mathbf{U}}_m}{{\hat D}_{jm}}}.\label{eq:InnerProductFinal} 
\end{equation}

Combining the discrete inner product result \eqref{eq:InnerProductFinal} with \eqref{eq:Hdot}, we determine the elemental values of $\dot H$ at a point $\left(\xi_{i},\eta_{j},\zeta_{k}\right)$  are 
\begin{equation}
\begin{split}
{{\dot {\mathbf{H}}}_{ijk}} &+ \left\{ {\left[ {{{\mathbf{F}}^*}\left( {1,{\eta _j},{\zeta _k}} \right)\frac{{{\ell _i}(1)}}{{{w_i}}} - {{\mathbf{F}}^*}\left( { - 1,{\eta _j},{\zeta _k}} \right)\frac{{{\ell _i}( - 1)}}{{{w_i}}}} \right] + \frac{1}{2}\sum\limits_{n=0}^N {\left[ {{{\left( {{{\widetilde {\mathbf{F}}}^1}} \right)}_{njk}} + {{\left( {{\mathcal{A}^1}} \right)}_{ijk}}{\mathbf{Q}_{njk}}} \right]{{\hat D}_{in}}} } \right\} \\&+ \left\{ {\left[ {{{\mathbf{F}}^*}\left( {{\xi _i},1,{\zeta _k}} \right)\frac{{{\ell _j}(1)}}{{{w_j}}} - {{\mathbf{F}}^*}\left( {{\xi _i}, - 1,{\zeta _k}} \right)\frac{{{\ell _j}( - 1)}}{{{w_j}}}} \right] + \frac{1}{2}\sum\limits_{n=0}^N {\left[ {{{\left( {{{\widetilde {\mathbf{F}}}^2}} \right)}_{ink}} + {{\left( {{\mathcal{A}^2}} \right)}_{ijk}}{\mathbf{Q}_{ink}}} \right]{{\hat D}_{jn}}} } \right\} \\&+ \left\{ {\left[ {{{\mathbf{F}}^*}\left( {{\xi _i},{\eta _j},1} \right)\frac{{{\ell _k}(1)}}{{{w_k}}} - {{\mathbf{F}}^*}\left( {{\xi _i},{\eta _j}, - 1} \right)\frac{{{\ell _k}( - 1)}}{{{w_k}}}} \right] + \frac{1}{2}\sum\limits_{n=0}^N {\left[ {{{\left( {{{\widetilde {\mathbf{F}}}^3}} \right)}_{ijn}} + {{\left( {{\mathcal{A}^3}} \right)}_{ijk}}{\mathbf{Q}_{ijn}}} \right]{{\hat D}_{kn}}} } \right\} \\&+ {G_{ijk}}{{\mathbf{Q}}_{ijk}} = 0\quad i,j,k=0,1,\ldots,N.
\end{split}
\label{eq:ImplementationHDot}
\end{equation}
We note that the structure of the skew-symmetric approximation clearly reflects its tensor product nature and is the same as the standard DGSEM  \cite{Kopriva:2009nx}, but now the interior flux terms are modified by adding terms of the form ${{\left( {{\mathcal{A}}} \right)}_{ijk}}{\mathbf{Q}_{njk}}$. Therefore,  the conversion to skew-symmetric form requires only including the additional terms and the factors of 1/2. Efficiency is gained by computing and storing those matrix-vector products line by line before the differentiation process, since for each node $(ijk)$ the coefficient matrix does not change.

Finally, we compute from $\dot{\mathbf{H}}$ (computed by (\ref{eq:ImplementationHDot})) and $\dot J$ (computed from the GCL, (\ref{eq:PointwiseJDot})) the time derivatives of the solution $\dot {\mathbf{Q}}$ or $\dot {J\mathbf{Q}}$, depending on whether the solution or the volume weighted solution is stored as the fundamental variable. Since at any nodal point $ijk$
\begin{equation}
{\dot {\mathbf{H}}_{ijk}} = \frac{1}{2}\left\{ {\frac{{\partial {J_{ijk}}{{\mathbf{Q}}_{ijk}}}}{{\partial \tau }} + {J_{ijk}}\frac{{\partial {{\mathbf{Q}}_{ijk}}}}{{\partial \tau }}} \right\} = {J_{ijk}}\frac{{\partial {{\mathbf{Q}}_{ijk}}}}{{\partial \tau }} + \frac{1}{2}\frac{{\partial J}}{{\partial \tau }}{{\mathbf{Q}}_{ijk}},\end{equation}
we can compute 
\begin{equation}
\dot{\mathbf{Q}}_{ijk} = \frac{{{{\dot {\mathbf{H}}}_{ijk}} - \frac{1}{2}{{\dot J}_{ijk}}{{\mathbf{Q}}_{ijk}}}}{{{J_{ijk}}}}\quad or\quad {\dot{\left( {J{\mathbf{Q}}} \right)}_{ijk}} = {{\dot {\mathbf{H}}} _{ijk}} + \frac{1}{2}{{\dot J}_{ijk}}{{\mathbf{Q}}_{ijk}}.
\end{equation}

With either form of the time derivative of the state vector, working backwards as in Sec. \ref{sec:GCLApprox} shows that the approximation remains constant state preserving. For example, when we multiply the left hand form by $J_{ijk}W_{ijk}$, sum over all points, and then replace $\dot{\mathbf{H}}$ and $\dot J$ by their approximations (\ref{eq:Hdot}) and (\ref{eq:DGGCL}), then
\[
\begin{split}
{\left( {J{\mathbf{\dot Q}},\boldsymbol\phi } \right)_N} &= \left( {\dot {\mathbf{H}} - \frac{1}{2}\dot J{\mathbf{Q}},\boldsymbol\phi } \right) \\&= -{\left. {\left\{ {{{\left( {{{ {\mathbf{F}} }^*} - \tilde{\mathbf{ F}}\cdot\hat n} \right)}^T}} \right\}\boldsymbol\phi } \right|_{\partial E,N}} - \frac{1}{2}\sum\limits_{i = 1}^3 {{{\left( {\frac{{\partial {{ {\mathbf{F}} }^i}}}{{\partial {\xi ^i}}},\boldsymbol\phi } \right)}_N}}  - \frac{1}{2}\sum\limits_{i = 1}^3 {{{\left( {\frac{{\partial {\mathbf{Q}}}}{{\partial {\xi ^i}}},\boldsymbol\phi } \right)}_N}}  + \frac{1}{2}{\left( {\sum\limits_{i = 1}^3 {\frac{{\partial {\mathbb{I}^N}\left( {{\mathcal{A}^i}} \right)}}{{\partial {\xi ^i}}}{\mathbf{Q}}} ,\boldsymbol\phi } \right)_N}.
\end{split}
\]
We now take the solution ${\mathbf{Q}} = \mathbf{c}=const$. When we do that, the boundary flux terms cancel due to consistency, the flux includes the constant state, and the derivative of the solution vanishes leaving
\[{\left( {J{{\dot {\mathbf{Q}}}},\boldsymbol\phi } \right)_N} =  - \frac{1}{2}\sum\limits_{i = 1}^3 {{{\left( {\frac{{\partial {\mathbb{I}^N}\left( {{\mathcal{A}^i}} \right){\mathbf{c}}}}{{\partial {\xi ^i}}},\boldsymbol\phi } \right)}_N}}  + \frac{1}{2}{\left( {\sum\limits_{i = 1}^3 {\frac{{\partial {\mathbb{I}^N}\left( {{\mathcal{A}^i}} \right)}}{{\partial {\xi ^i}}}{\mathbf{c}}} ,\boldsymbol\phi } \right)_N} = 0.\]
Choosing $\boldsymbol\phi$ to be the Lagrange basis polynomials shows that $\dot {\mathbf{Q}} =0$ at each nodal point.

The  equations for the time derivatives of the state vector and Jacobian are then integrated with an explicit time integrator, e.g. Runge-Kutta or Adams-Bashforth methods. 

\section{Numerical Results}\label{NumericalResultsSection}

We provide four examples that combine the skew-symmetric DGSEM-ALE spatial discretization with the low-storage third order Runge-Kutta scheme of Williamson \cite{Williamson:1980:JCP80}. The numerical tests are selected to demonstrate the theoretical properties of the numerical scheme derived in Sec. \ref{PropertiesOfStable}. The first demonstrates the high-order convergence in space and full time accuracy of the approximation on a moving mesh. The second supports the stability proven in Sec. \ref{sec:ProvableStability}. We show a case where the energy growth for the skew-symmetric approximation remains bounded in time while that of the standard DGSEM-ALE approximation does not. Next, we demonstrate that the skew-symmetric approximation remains conservative as shown in Sec. \ref{ConservationSection}. Finally, we numerically verify free-stream preservation, the result proved in Sec. \ref{FreeStreamSection}.

For each of the numerical tests in this paper we consider the symmetric form of the moving mesh wave equation written as a conservation law

\begin{equation}\label{waveeqn} 
\begin{aligned}
\left[
\begin{array}{c}
  p \\
  u \\
  v \\
  w
\end{array}
\right]_\tau &+ 
 \left[\begin{array}{cccc}
  -x_\tau & c & 0 & 0 \\
  c &  -x_\tau  & 0 & 0 \\
  0 &  0  & -x_\tau & 0 \\
  0 & 0 & 0 & -x_\tau
\end{array}
\right]
 \left[\begin{array}{c}
  p \\
  u \\
  v \\ 
  w
\end{array}
\right]_x + 
 \left[\begin{array}{cccc}
  -y_\tau &  0  & c & 0 \\
  0 &  -y_\tau  & 0 & 0 \\
  c &  0  & -y_\tau & 0 \\
  0 & 0 & 0 & -y_\tau
\end{array}
\right]
 \left[\begin{array}{c}
  p \\
  u \\
  v \\
  w
\end{array}
\right]_y \\
&\qquad\qquad + \left[\begin{array}{cccc}
  -z_\tau &  0  & 0 & c \\
  0 &  -z_\tau  & 0 & 0 \\
  0 &  0  & -z_\tau & 0 \\
  c & 0 & 0 & -z_\tau
\end{array}
\right]
 \left[\begin{array}{c}
  p \\
  u \\
  v \\
  w
\end{array}
\right]_z = 0,
\end{aligned}
\end{equation}
where $\mathbf{x}_\tau = (x_\tau,y_\tau,z_\tau)$ is the mesh velocity and $c$ is a constant wave speed. For each numerical test we take the wave speed $c=1$.

\begin{figure}[!ht]
\begin{center}
  \includegraphics[scale=0.5]{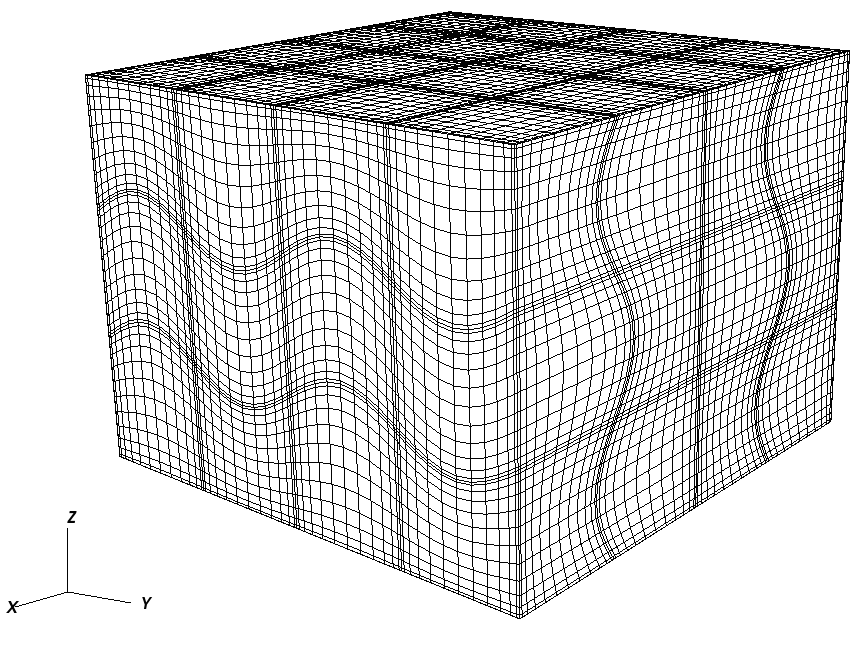}
\caption{The curvilinear mesh used for verification of convergence, stability, conservation, and free-stream preservation.}
\label{fig:PeriodicMesh}
\end{center}
\end{figure}

The skew-symmetric DGSEM-ALE approximation is general and operates on unstructured meshes. For convenience, however, we will use a structured curvilinear mesh, where it is straightforward to set periodic boundary conditions for testing conservation, to demonstrate the theoretical attributes of the approximation. We consider the domain $\Omega = [-2,2]\times[-2,2]\times[0,3]$ divided into 48 elements. The element lengths are uniform in each Cartesian direction. We then create a periodic curvilinear mesh by replacing flat planes, e.g. $x = 1$, with a sinusoidal plane. An example of the curvilinear mesh with $N=12$ in each spatial direction is shown in Fig. \ref{fig:PeriodicMesh}.

Finally, we describe the mesh motion used for each of the test problems. We prescribe a periodic motion to the corner nodes of a hexahedral element that initially lies on the flat plane $y=0$. In this way we have an analytical representation of the mesh motion and mesh velocities of the moving corner nodes. With this information we can update the curvilinear hexahedral element geometry while maintaining high-order accuracy \cite{Minoli:2010rt}. The corner node motion is given by
\begin{equation}\label{meshMotion}
\begin{aligned}
x &= {x}_* -\frac{1}{4}\sin(2\pi \tau) \\ 
y &= {y}_* +\frac{1}{4}\sin(2\pi \tau) \\ 
z &= {z}_* +\frac{1}{4}\sin(2\pi \tau) \\ 
\end{aligned}
\end{equation}
where $x_*$, $y_*$, and $z_*$ are the initial positions for a corner node. It is straightforward to obtain analytical expression for the mesh velocities $\mathbf{x}_\tau$ from \eqref{meshMotion}.

\subsection{Convergence}

For the numerical convergence study we choose initial and boundary conditions so that the solution is a Gaussian plane wave
\begin{equation}\label{planeWave}
 \left[
\begin{array}{c}
  p \\
  u \\
  v \\
  w
\end{array}
\right] = 
 \left[\begin{array}{c}
  1 \\
  \frac{k_x}{c} \\
  \frac{k_y}{c} \\
  \frac{k_z}{c}
\end{array}
\right]e^{-\frac{(k_x(x-x_0)^2+k_y(y-y_0)^2+k_z(z-z_0)^2-ct)^2}{d^2}},
\end{equation}
with the wavevector $\vec{k}$ normalized to satisfy $k_x^2+k_y^2 + k_z^2= 1$, $d = 1.0$, and vary $x_0$, $y_0$, $z_0$ to adjust the initial position of the wave.

\begin{figure}[!ht]
\begin{center}
  \includegraphics[scale=0.75]{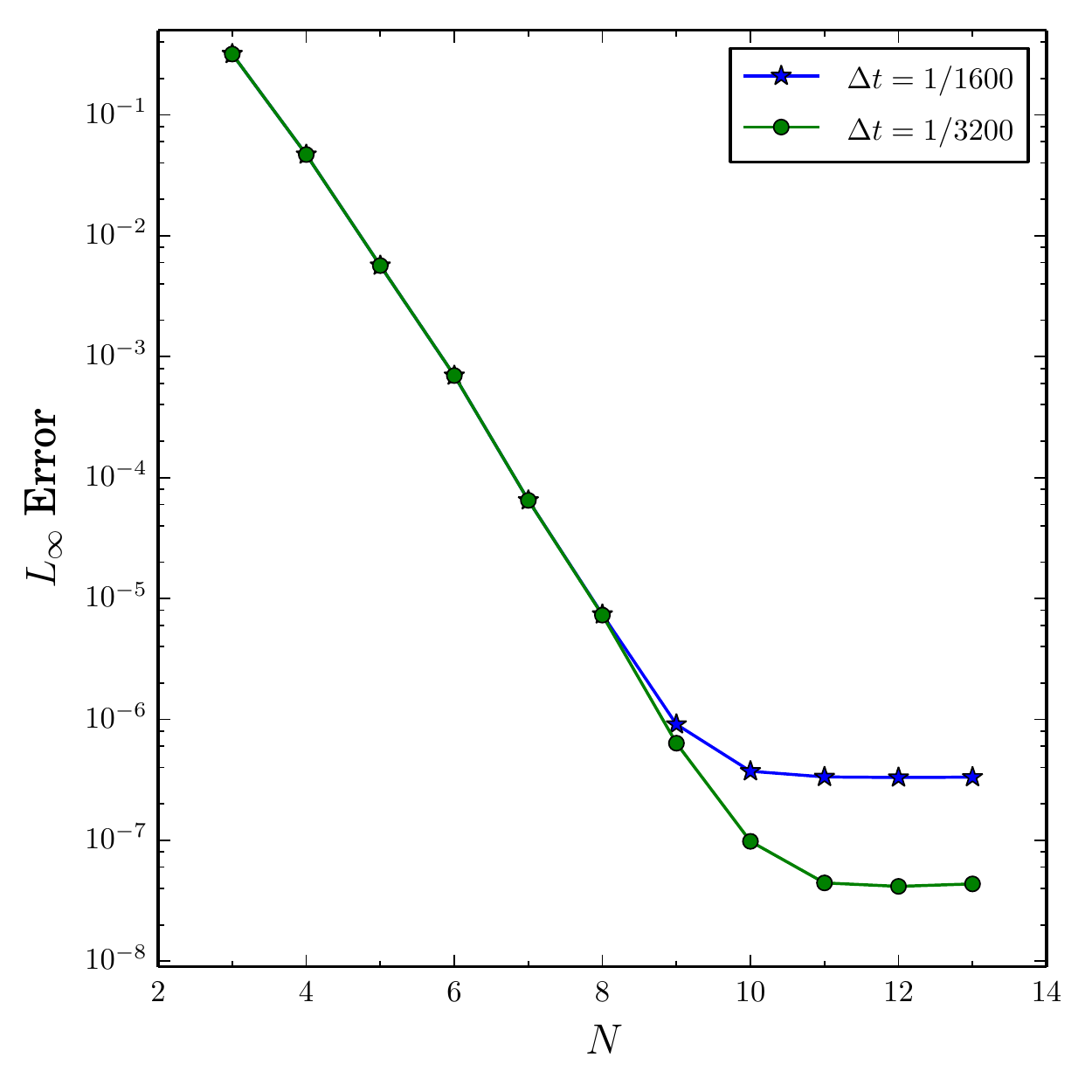}
\caption{Semi-log plot visualizes the spectral convergence for the skew-symmetric DGSEM-ALE approximation. When $\Delta t$ is halved, the error for large $N$ drops by a factor of eight, which demonstrates the third order temporal accuracy.}
\label{fig:Convergence}
\end{center}
\end{figure}

We obtain spectral accuracy in space and design accuracy in time for the skew-symmetric DGSEM-ALE. For the example, we take the initial position to be $x_0=-1.0$ and $y_0=z_0=0.0$. We integrate the solution up to the final time $T=4.0$. Fig. \ref{fig:Convergence} shows exponential convergence in space until $N = 10$, where the error becomes dominated by time integrator errors. Note that when the value of $\Delta t$ is halved the error in the approximation is reduced by a factor of 8, as expected for the third order time integration technique.

\subsection{Stability}

\begin{figure}[!ht]
\begin{center}
  \includegraphics[scale=0.75]{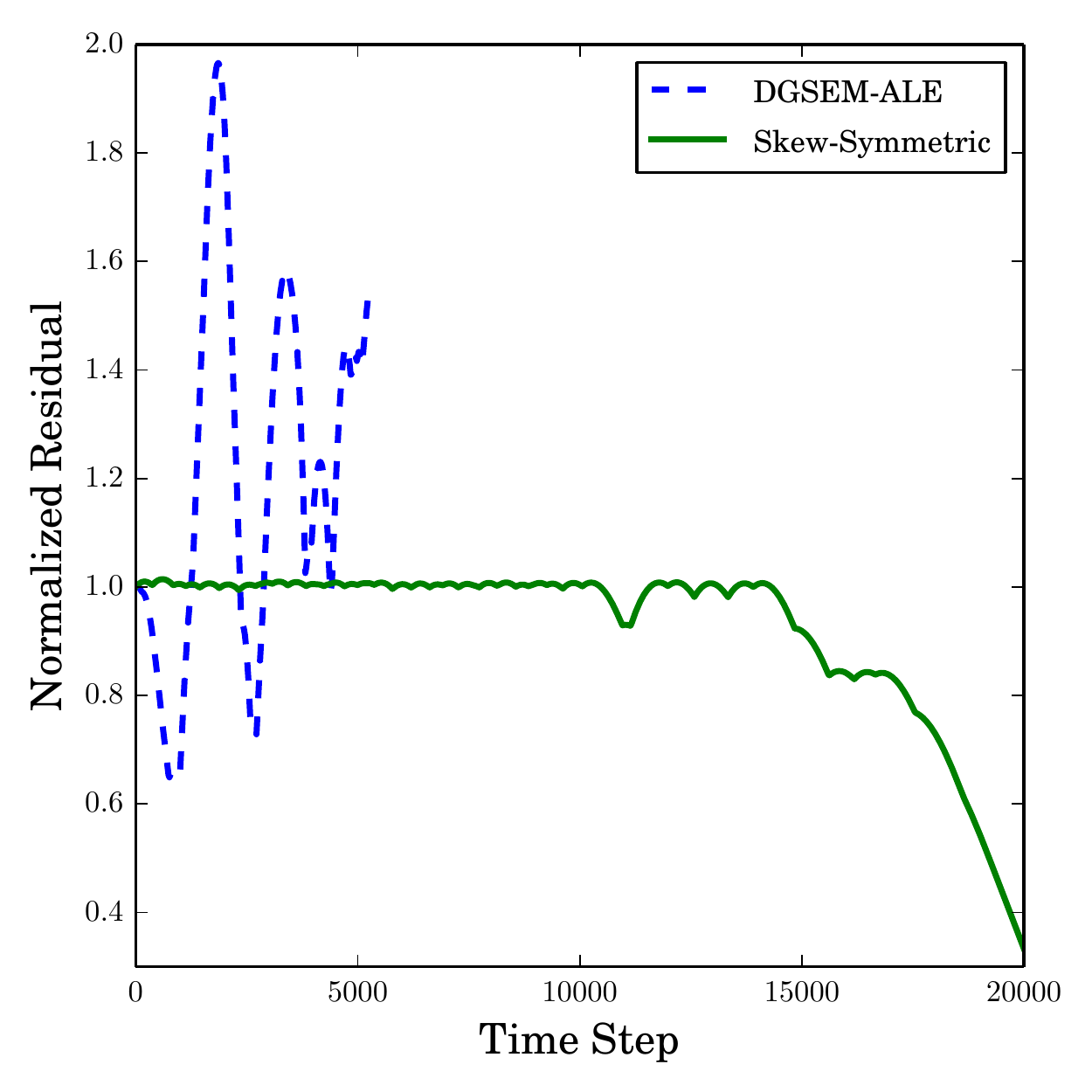}
\caption{The normalized maximum residual as a function of time demonstrates the stability of the skew-symmetric formulation and the weak instability of the classical DGSEM approximation for the moving mesh in Fig. \ref{fig:PeriodicMesh}.}
\label{fig:Stability}
\end{center}
\end{figure}

We next demonstrate the stability of the skew-symmetric formulation. We consider an identical plane wave configuration \eqref{planeWave} to show instability of a standard DGSEM-ALE for this problem. The computation is run for twenty thousand time steps with $T=6$ as the final time. We compute the maximum residual, $(JQ)_\tau$, normalized to the initial value, as a function of time for the central numerical flux with polynomial order $N=4$ and show it in Fig. \ref{fig:Stability}. We observe that the residual for the skew-symmetric formulation remains bounded and decreases as the wave moves out of the domain, whereas the standard DGSEM-ALE is weakly unstable, blowing up after about six thousand time steps. These results illustrate the proposition as presented in Sec. \ref{sec:ProvableStability} that the approximation is stable.

\subsection{Conservation}

We next show an example of the result of Sec. \ref{ConservationSection} that the skew-symmetric DGSEM-ALE remains conservative despite the addition of the non-conservative form of the equation. To do so we integrate the wave equation \eqref{waveeqn} with the initial condition
\begin{equation}
p(x,y,z,0) = \exp\left[-\ln(2)\left(\frac{x^2+y^2+z^2}{2.3^3}\right)\right],\;u(x,y,z,0) = v(x,y,z,0) = w(x,y,z,0) = 0,
\end{equation}
and periodic boundary conditions. We compute the solution up to a final time $T=1$ with a fixed time step $\Delta t = 1/1000$. To test conservation we compute the total amount of each of the conserved quantities over the entire domain, e.g.
\begin{equation}
p^{tot} = \sum_{m=1}^{N_{el}}\sum_{i,j,k=0}^NJ^m_{ijk}p^m_{ijk}W_iW_jW_k.
\end{equation}
Tables \ref{tab:ConservationUpwind} and \ref{tab:ConservationCentral} present the $\mathbb{L}^{\infty}$ error in the total amount of each conserved variable compared with the initial value for the upwind and central numerical fluxes respectively. All the computed errors are on the order of double precision roundoff for the polynomial orders $N=3$ and $N=4$. These results support the result from Sec. \ref{ConservationSection} that the moving mesh skew-symmetric approximation is globally conservative.
\begin{table}[!ht]
\begin{center}
  \begin{tabular}{lcc}
  \toprule
  $\mathbb{L}^{\infty}$ Error & $N=3$ & $N=4$ \\
 \midrule
 $\|p^{tot} - p^{tot}_0\|_{\infty}$ & $7.00 \times 10^{-15}$ &  $1.42 \times 10^{-14}$  \\
 $\|u^{tot} - u^{tot}_0\|_{\infty}$ & $5.19 \times 10^{-15}$ &  $5.76 \times 10^{-15}$  \\
 $\|v^{tot} - v^{tot}_0\|_{\infty}$ & $1.17 \times 10^{-15}$ &  $2.55 \times 10^{-15}$  \\
 $\|w^{tot} - w^{tot}_0\|_{\infty}$ & $4.60 \times 10^{-15}$ &  $8.71 \times 10^{-16}$  \\
 \bottomrule
  \end{tabular}
  \caption{$\mathbb{L}^{\infty}$ error of the total amount of each conserved quantity demonstrates the global conservative property of the approximation with the upwind numerical flux.}
  \label{tab:ConservationUpwind}
\end{center}
\end{table}

\begin{table}[!ht]
\begin{center}
  \begin{tabular}{lcc}
  \toprule
  $\mathbb{L}^{\infty}$ Error & $N=3$ & $N=4$ \\
 \midrule
 $\|p^{tot} - p^{tot}_0\|_{\infty}$ & $1.40 \times 10^{-14}$ &  $1.41 \times 10^{-14}$  \\
 $\|u^{tot} - u^{tot}_0\|_{\infty}$ & $5.29 \times 10^{-15}$ &  $5.36 \times 10^{-15}$  \\
 $\|v^{tot} - v^{tot}_0\|_{\infty}$ & $1.02 \times 10^{-15}$ &  $2.54 \times 10^{-15}$  \\
 $\|w^{tot} - w^{tot}_0\|_{\infty}$ & $4.32 \times 10^{-15}$ &  $9.71 \times 10^{-16}$  \\
 \bottomrule
  \end{tabular}
  \caption{$\mathbb{L}^{\infty}$ error of the total amount of each conserved quantity demonstrates the global conservative property of the approximation with the central numerical flux.}
  \label{tab:ConservationCentral}
\end{center}
\end{table}

\subsection{Free-Stream Preservation}

Finally, to show that the skew-symmetric DGSEM-ALE preserves a constant solution when a mesh moves, we consider a uniform solution in time
\begin{equation}\label{contSolution}
\mathbf{q} =  \left[
\begin{array}{c}
  p \\
  u \\
  v \\
  w
\end{array}
\right] =  \left[
\begin{array}{c}
  \pi \\
  \pi \\
  \pi \\
  \pi
\end{array}
\right].
\end{equation}
We compute the solution on the mesh shown in Fig. \ref{fig:PeriodicMesh} and prescribe the mesh motion by \eqref{meshMotion}.

The error for the test of free-stream preservation was calculated using the maximum norm of the computed solution against the exact constant solution at $T = 2.0$, which corresponds to a complete period in the oscillation of the vertical plane $y=0$. We fix the time step to be $\Delta t = 1/1000$. Table \ref{Tab:FreeStream} shows the computed error for double precision computations for two polynomial orders, $N=3$ and $N=4$. We found in the computations that the constant state is preserved for either the upwind, $\lambda = 1$, or central, $\lambda =0$, numerical flux.
\begin{table}[!ht]
\begin{center}
  \begin{tabular}{lcc}
  \toprule
   $\mathbb{L}^{\infty}$ Error & $N=3$ & $N=4$ \\
 \midrule
 Upwind $(\lambda=1)$ & $3.97 \times 10^{-13}$ &  $4.16 \times 10^{-13}$  \\
 Central $(\lambda=0)$ & $3.97 \times 10^{-13}$ &  $4.16 \times 10^{-13}$  \\
 \bottomrule
  \end{tabular}
  \caption{$\mathbb{L}^{\infty}$ error of a constant solution to the wave equation computed with the skew-symmetric DGSEM on a moving mesh.}
  \label{Tab:FreeStream}
\end{center}
\end{table}

\section{Conclusion}\label{ConclusionSection}

We have derived a new, provably energy stable DGSEM formulation for moving elements. The motion of the domain was handled using a time dependent mapping in the ALE framework, where systems of conservation laws on a moving domain remain systems of conservation laws on a static reference domain. To create an energy stable approximation we used the skew-symmetric form of the governing equations. We proved that the approximation is stable, constant state preserving, and remains globally conservative despite the addition of the non-conservative form of the equation. Numerical studies support the theoretical results.



\bibliographystyle{plain}
\bibliography{dakBib.bib}

\end{document}